\documentclass[12pt, reqno]{amsart}
\usepackage{amsmath, amstext, amsbsy, amssymb, amscd}

\newtheorem{theorem}{Theorem}[section]
\newtheorem{Definition-Lemma}[theorem]{Definition-Lemma}
\newtheorem{lemma}[theorem]{Lemma}
\newtheorem{proposition}[theorem]{Proposition}
\newtheorem{corollary}[theorem]{Corollary}

\theoremstyle{definition}
\newtheorem{definition}[theorem]{Definition}

\newtheorem{example}[theorem]{Example}
\newtheorem{conjecture}[theorem]{Conjecture}

\theoremstyle{remark}
\newtheorem{remark}[theorem]{Remark}

\numberwithin{equation}{section}

\newcommand{\ale}{\C^2//\G}
\newcommand{\be}{\beta}

\newcommand{\C}{ \mathbb C }

\newcommand{\g}{{\gamma}}
\newcommand{\G}{{\Gamma}}
\newcommand{\Gn}{{\G_n}}

\newcommand{\la}{\lambda}

\newcommand{\N}{ \mathbb N }

\newcommand{\Rgn}{R(\Gn)}
\newcommand{\Xn}{X^{[n]} }
\newcommand{\Z}{ \mathbb Z }

{\vskip-\lastskip\medskip
  \noindent
  {\em #1.}\enspace
  }%
{\qed\par\medskip
  }

\begin{document}

\title[The FH-ring of wreath products and
Hilbert schemes]
{The Farahat-Higman ring of wreath products and
Hilbert schemes}

\author[Weiqiang Wang]{Weiqiang Wang}
\address{Department of Mathematics, University of Virginia,
Charlottesville, VA 22904} \email{ww9c@virginia.edu}

\begin{abstract}

We study the structure constants of the class algebra $R_\Z(\Gn)$
of the wreath products $\Gn$ associated to an arbitrary finite
group $\G$ with respect to the basis of conjugacy classes. We show
that a suitable filtration on $R_\Z(\Gn)$ gives rise to the graded
ring $\mathcal G_\G(n)$ with non-negative integer structure
constants independent of $n$ (some of which are computed), which
are then encoded in a Farahat-Higman ring $\mathcal G_\G$. The
real conjugacy classes of $\G$ come to play a distinguished role,
and is treated in detail in the case when $\G$ is a subgroup of
$SL_2(\C)$. The above results provide new insight to the
cohomology rings of Hilbert schemes of points on a
quasi-projective surface $X$.
\end{abstract}

\maketitle

\date{}


\section{Introduction}

\subsection{}
The wreath products $\Gn  =\G^n \rtimes S_n$ associated to a
finite group $\G$ are natural generalizations of the symmetric
groups $S_n$ (cf. \cite{Zel, Mac}). Connections of the wreath
products with Hilbert schemes of points on surfaces and with
Nakajima quiver varieties \cite{Na1, Na3} were first pointed out
by the author (cf. \cite{Wa3} for an overview). For example, it
has since been expected that all geometric invariants on the
Hilbert scheme $(\ale)^{[n]}$ of $n$ points on the minimal
resolution $\ale$ associated to $\G\leq SL_2(\C)$ can be entirely
described using the wreath product $\Gn$. The interrelations among
wreath products, Hilbert schemes, and Nakajima quiver varieties
have been one of the main motivations for us to study various
aspects of wreath products, cf. \cite{Wa1, FJW, Wa2, WaZ}. Also
see \cite{EG, Wa4} for more recent development.

The first aim of this paper is to develop a new approach of
studying the class algebras of the wreath products associated to
an arbitrary finite group $\G$. We establish various properties of
the structure constants of the class algebras $R_\Z(\Gn)$ with
respect to the basis of conjugacy classes. This allows us to
introduce a ring $\mathcal G_\G$ which encodes all the structures
of the graded algebras associated to a natural filtration on
$R_\Z(\Gn)$ for all $n$. We call $\mathcal G_\G$ the
Farahat-Higman ring, or the FH-ring, since our results specialize
to the ones of Farahat-Higman \cite{FH} for the symmetric groups
(see Macdonald \cite{Mac} for an elegant presentation and
improvement). While we are much inspired by the approach of
\cite{FH}, we need to generalize various concepts from the
symmetric groups and develop new delicate combinatorial analysis
in order to treat the extra complications coming from the presence
of the group $\G$. The appearance of new phenomena which are not
observable in the symmetric group case will be crucial for the
subsequent geometric applications.

The second aim of this paper is to point out the geometric
implications of these results on wreath products. The FH-ring,
when $\G$ is a finite subgroup of $SL_2(\C)$, is isomorphic to the
cohomology rings of Hilbert schemes of points on the minimal
resolutions $\ale$ of a simple singularity $\C^2/\G$. In this way,
our concrete results on the class algebras of $\Gn$ provide new
insights into the cohomology rings of the corresponding Hilbert
schemes. Perhaps no less significantly, the results on the
Farahat-Higman ring etc enable us to predict some remarkable
general structures of the cohomology rings of Hilbert schemes
$\Xn$, most notably for an arbitrary quasi-projective surface $X$.
Let us explain in more detail.

\subsection{}

The conjugacy classes of the wreath product $\Gn$ are determined
by their types, which are in one-to-one correspondence with the
set $\mathcal P_n(\G_*)$ of partition-valued functions on $\G_*$
such that the total number of boxes in the corresponding Young
diagrams is $n$. Here and further $\G_*$ denotes the set of
conjugacy classes of $\G$. The union $\G_\infty  =\cup_{n \ge 1}
\G_n$ carries a group structure by the natural embedding of $\Gn$
in $\G_{n+1}$. We introduce a notion of {\em modified types}, (cf.
\cite{Mac}, pp. 131 for the symmetric group case), so that the set
$\mathcal P (\G_*)$ of modified types (which are all
partition-valued functions on $\G_*$) parametrizes naturally the
conjugacy classes of $\G_\infty$. By assigning a certain
non-negative integer to an element or its modified type, we define
a function $\| \cdot \|$ on $\G_\infty$ or alternatively on
$\mathcal P (\G_*)$. We introduce a notion of {\em reduced
expression} for an element in $\G_\infty$, and further identify
$\|x\|$ with the length of a reduced expression for $x \in
\G_\infty$. The function $\| \cdot \|$ is sub-multiplicative and
thus defines a filtration on the class algebra $R_\Z (\Gn)$ for
each $n$.

Let us denote by $K_\la (n)$ the sum of elements in the conjugacy
class in $\Gn$ of modified type $\la$, and write the
multiplication $K_\la (n) \, K_\mu(n) =\sum_\nu a_{\la\mu}^\nu (n)
\, K_\nu(n)$ where $a_{\la\mu}^\nu (n) \in \Z_+$, and
$a_{\la\mu}^\nu (n) =0$ unless $\| \nu \| \le \| \la \|+ \| \mu
\|$. We prove that the structure constants $a_{\la\mu}^\nu (n)$
are polynomials in $n$ (whose degree is also described
explicitly); and when $\| \nu \| = \| \la \|+ \| \mu \|$,
$a_{\la\mu}^\nu (n)$ is a constant $a_{\la\mu}^\nu \in \Z_+$
independent of $n$. This allows us to introduce a commutative ring
$\mathcal F_\G$ with a basis $(K_\la)_{\la \in \mathcal P(\G_*)}$
and a $\Z_+$ filtration, whose associated graded ring $\mathcal
G_\G$ (called the {\em Farahat-Higman ring}) has a multiplication
$K_\la  \, K_\mu =\sum_{\| \nu \| = \| \la \|+ \| \mu \|}
a_{\la\mu}^\nu  \; K_\nu$. Similarly, the filtration on
$R_\Z(\Gn)$ also gives rise to a graded ring, denoted by $\mathcal
G_\G (n)$, whose structure constants are non-negative integral and
independent of $n$. The ring $\mathcal G_\G$ determines a family
of rings $\{\mathcal G_\G (n)\}_{n\ge 1}$ and vice versa. When
$\G$ is trivial, we write $\mathcal G_\G (n)$ as $\mathcal G (n)$.

We further calculate explicitly these structure constants
$a_{\la\mu}^\nu$ when $\mu$ is of single-cycle modified type. The
calculation requires a rather careful combinatorial analysis of
the convolution products for $\Gn$, which is made possible by the
rigid constraint $\| \nu \| = \| \la \|+ \| \mu \|$. Let us
mention below five interesting consequences of this computation.
First, it provides us an algorithm for computing the general
structure constants. Secondly, the ring $\mathbb Q \otimes_\Z
\mathcal G_\G$ is shown to be a polynomial ring generated by
$K_\mu$, where $\mu$ runs over all possible single-cycle modified
types. Thirdly, we obtain a simple set of ring generators for the
ring $\mathcal G_\G(n)$. Fourthly, we show that the ring $\mathcal
G_\G (n)$ for any $\G$ contains the ring $\mathcal G (n)$ as a
natural quotient. Fifthly, we observe a remarkable fact that the
new structure constants of the ring $\C \otimes_\Z \mathcal G_\G$
(over $\C$!) with respect to a suitably rescaled basis depend no
longer on the group $\G$, but only on two integers: the
cardinalities of $\G_*$ and its subset $\G_*^{re}$ of nontrivial
real conjugacy classes.

\subsection{}

It is well known that the set of finite subgroups $\G$ of
$SL_2(\C)$ corresponds bijectively to the Dynkin diagrams of ADE
types. The dual McKay correspondence of Ito-Reid \cite{IR} (also
cf. Brylinski \cite{Bry}) provides further a canonical bijection
between the nontrivial conjugacy classes of $\G$ and vertices of
the associated Dynkin diagram. For every finite subgroup $\G$ of
$SL_2(\C)$, we determine the set $\G_*^{re}$ and further observe
that $\G_*^{re}$ affords an elegant interpretation as the
fixed-point set of a distinguished Dynkin diagram automorphism,
via the dual McKay correspondence. This automorphism in turn has
simple interpretations in surface geometry as well as in Lie
theory. To our best knowledge, the significance of the real
conjugacy classes was not noted before in the literature on the
McKay correspondence. We remark that there has been direct
connection between wreath products and McKay correspondence
\cite{FJW, Wa2}. However a direct interplay between McKay and dual
McKay correspondences has yet to be discovered.

Note that $\Gn$ acts on $\C^{2n}$ naturally, and the orbifold
$\C^{2n} /\Gn $ can be regarded as a generalization of the simple
singularity $\C^2/\G$. The author constructed a crepant resolution
of singularities $ (\ale)^{[n]} \rightarrow \C^{2n} /\Gn $ and
another crepant resolution of singularities of $\C^{2n} /\Gn $
which can be identified with a Nakajima variety (cf. \cite{Wa1,
Wa2, Wa3}). A combination of a theorem of Etingof-Ginzburg
\cite{EG} and results of Nakajima \cite{Na1} etc\footnote{A quiver
variety $\mathbb M_{\G,n}$ is used in \cite{EG} in place of the
Hilbert scheme $(\ale)^{[n]}$. However, according to Nakajima
(unpublished) and Kuznetsov \cite{Kuz}, $(\ale)^{[n]}$ affords a
same quiver variety description with different stability
conditions. By a theorem in \cite{Na1}, $(\ale)^{[n]}$ and
$\mathbb M_{\G,n}$ are diffeomorphic, and thus have isomorphic
cohomology rings.} says that there is a ring isomorphism between
the cohomology ring $H^*( (\ale)^{[n]})$ with $\C$-coefficient and
the graded ring associated to the orbifold filtration on $\C
\otimes_\Z R_\Z(\Gn)$ defined by the shift numbers \cite{Zas}
(also cf. \cite{IR}). These shift numbers for $\Gn$ were computed
earlier in \cite{WaZ}, and they coincides with the degree $\|
\cdot \|$ introduced in this paper for $\G \leq SL_2(\C)$. Thus,
there is a ring isomorphism between $H^*( (\ale)^{[n]})$ and
$\C\otimes_\Z \mathcal G_\G(n)$. Note that the ring $\C\otimes_\Z
\mathcal G_\G(n)$ for $\G \leq SL_2(\C)$ can alternatively be
interpreted as the Chen-Ruan cohomology ring of $\C^{2n}/\Gn$, and
the above ring isomorphism supports conjectures of Ruan (cf.
\cite{Ru}). In the important special case $X=\C^2$, i.e. $\G$ is
trivial, the ring isomorphism was due to Lehn-Sorger \cite{LS1}
and independently Vasserot \cite{Vas} (but the connections to
\cite{FH} was not noticed). In this way, our results on $\mathcal
G_\G(n)$ and $\mathcal G_\G$ for $\G \leq SL_2(\C)$ can be
regarded as new results on the ring $H^*( (\ale)^{[n]})$.

\subsection{}

In recent years, there has been much activity in understanding the
cohomology rings of Hilbert schemes $\Xn$ of $n$ points on a
projective surface $X$. However, when $X$ is quasi-projective (the
usage of the terminology {\em quasi-projective} in this paper
excludes {\em projective}), the connections with vertex operators
become obscure, and little is known about the cohomology rings of
$\Xn$ with the exception of the affine plane. Our guiding
principle here (cf. \cite{Wa3}) is that the study of the class
algebras of $\Gn$ associated to an arbitrary finite group $\G$
(resp. the graded ring associated to some suitable filtration)
should mirror the study of the cohomology rings of $\Xn$ of points
on an arbitrary surface $X$ which is projective (resp.
quasi-projective). From this view, the results on wreath products
in this paper should reflect new finer structures of the
cohomology rings of $\Xn$ for $X$ quasi-projective. This will be
indicated at the end of this paper, and our main conjecture says
that the cohomology rings of $\Xn$ for a quasi-projective surface
$X$ and all $n$ are governed by a Farahat-Higman type ring
associated to $X$. We hope that understanding the cohomology ring
of $\Xn$ for a quasi-projective surface $X$ may shed light on the
same problem for Nakajima quiver varieties, which is typically
noncompact. In another direction, we will address elsewhere
certain geometric and algebraic deformations of the current work
and \cite{Wa4}.

The paper is organized as follows. In Sect.~\ref{sec:FHring} we
set up the notations for wreath products, and then establish
various properties of the class algebras of the wreath products
which lead to the notion of a Farahat-Higman ring $\mathcal G_\G$.
In Sect.~\ref{sec:computation} we compute certain structure
constants for the FH-ring $\mathcal G_\G$ and describe various
implications of this computation. We further describe the
connections with the dual McKay correspondence. In
Sect.~\ref{sec:scheme} we present the connections with Hilbert
schemes, and describe some predictions on the cohomology rings of
$\Xn$ when $X$ is quasi-projective.

{\bf Acknowledgement.} I am very grateful to Zhenbo Qin for his
insightful collaboration on closely related projects as well as
numerous stimulating discussions on Hilbert schemes. I thank Igor
Dolgachev for helpful discussions and answering my questions about
surfaces (which in particular include the minimal resolutions in
Remark~\ref{rem:autom}) and his interest in the present work. This
work is partially supported by an NSF grant. I also thank MSRI for
its hospitality and providing an excellent working atmosphere in
Spring 2002, where this work was carried out.

\noindent{\em Note Added.} Conjecture~\ref{conj:indep} has been
established in some cases in \cite{LQW4}.

\section{The Farahat-Higman ring of wreath products}
\label{sec:FHring}
\subsection{Basics of the wreath products}

Let $\G$ be a finite group, and $\G_*$ be the set of
conjugacy classes of $\G$. The order of the centralizer
of an element in a conjugacy class $c\in\G_*$ is denoted
by $\zeta_c$. We will denote the identity conjugacy
class in $\G$ by $c^0$ and the identity of $\G$
by $1$. Denote by $\bar{c}$ the conjugacy
class $\{ x \mid x^{-1} \in c \}$. We will say that two elements
of $\G$ lie in {\em opposite} conjugacy classes, if
one belongs to some conjugacy class
$c \in \G_* \backslash c^0$ while the other belongs
to $\bar{c}$ (we allow that
$\bar{c}=c$). We shall denote by $\langle x \rangle$
the conjugacy class of an element $x$ in $\G$.

The symmetric group $S_n$ acts on the product group
$\Gamma^n= \Gamma \times \cdots
\times \Gamma$ by permutations: $\sigma
(g_1, \cdots, g_n)
  = (g_{\sigma^{ -1} (1)}, \cdots, g_{\sigma^{ -1} (n)}).
$ The wreath product of $\Gamma$ with $S_n$ is defined to be the
semidirect product
$$
 \Gamma_n = \{(g, \sigma) | g=(g_1, \cdots, g_n)\in {\Gamma}^n,
\sigma\in S_n \}
$$
 with the multiplication
\begin{eqnarray*}  \label{eq:twistprod}
 (g, \sigma)\cdot (h, \tau)=(g \, {\sigma} (h), \sigma \tau).
\end{eqnarray*}
The $i$-th factor subgroup of the product group $\G^n$ will be
denoted by $\G^{(i)}$. We have $\G^{(i)} \leq \G^n \leq \Gn$
and a short exact sequence of groups
$$
 1 \longrightarrow \G^n \longrightarrow \Gn
{\longrightarrow} S_n
 \longrightarrow 1.$$

The space $R_\Z(\G)$ of class functions of a finite group $\G$
(which is often called the {\em class algebra} of $\G$)
is closed under
the group multiplication (which is often referred to as
the {\em convolution product}). In this way, $R_\Z(\G)$
is also identified with the center of the group algebra $\Z[\G]$.
We will denote by $\Rgn$ (resp. $R_\Z(\Gn)$)
the {\em class algebra}
of complex (resp. integer) class functions on $\Gn$
endowed with the convolution product.

We denote by $|\la|=\la_1+\cdots+\la_\ell$ and the length $\ell(\la)
=\ell$, for a partition $\la=(\la_1, \la_2, \cdots, \la_\ell)$, where
$\la_1 \geq \dots \geq \la_\ell \geq 1$. We will also make use of
another notation for partitions:
$ \la=(1^{m_1(\la)}2^{m_2(\la)}\cdots) ,$
where $m_i(\la)$ is
the number of parts in $\la$ equal to $i$. Given two
partitions $\la, \mu$, we denote by $\la \supset \mu$
if $\mu$ is consisted of some parts of $\la$. In particular,
we have the empty partition $\emptyset \subset \la$
for any $\la$. For $\la \supset \mu$, we denote
by $\la -\mu$ the partition obtained by removing the parts
of $\mu$ from $\la$. We denote by $\la \cup \mu$
the partition whose parts are those of $\la$ and $\mu$,
arranged in descending order. We denote by $\la \geq \mu$
if $\la$ dominates $\mu$.

For a finite set $Y$ and $\rho=(\rho(x))_{x \in Y}$ a family of
partitions indexed by $Y$, the {\em degree} of $\rho$ is
$\|\rho\|=\sum_{x \in Y}|\rho(x)|$ and the {\em length} of $\rho$
is $\ell (\rho) = \sum_{c \in Y} \ell(\rho(c))$. It is convenient
to regard $\rho=(\rho(x))_{x \in Y}$ as a partition-valued
function on $Y$. We denote by ${\mathcal P}(Y)$ the set of all
partition-valued functions on $Y$ and put ${\mathcal P}_n (Y)=
\{\rho \in {\mathcal P} (Y) \mid \|\rho\| =n \}$. Given $\rho,
\sigma \in {\mathcal P}_n (Y)$, we extend naturally the
definitions from partitions to ${\mathcal P}_n (Y)$ by doing so
pointwise on $Y$ to define $\rho \supset \sigma$, $\rho - \sigma$
and $\rho \cup \sigma$, etc.

The conjugacy classes of $\Gn$ can be described in the following
way (cf. \cite{Mac}). Let $x=(g, \sigma )\in {\Gamma}_n$, where
$g=(g_1, \ldots , g_n) \in {\Gamma}^n,$ $ \sigma \in S_n$. The
permutation $\sigma $ is written as a product of disjoint cycles.
For each such cycle $y=(i_1 i_2 \cdots i_k)$, the element
$p_y =g_{i_k}
g_{i_{k -1}} \cdots g_{i_1} \in \Gamma$ is determined up to
conjugacy in $\Gamma$ by $g$ and $y$, and will be called the {\em
cycle-product} of $x$ corresponding to the cycle $y$. For any
conjugacy class $c \in \G_*$ and each integer $i \geq 1$, the number of
$i$-cycles in $\sigma$ whose cycle-product lies in $c$ will be
denoted by $m_i(\rho(c))$, and $\rho (c)$ denotes the partition
$(1^{m_1 (c)} 2^{m_2 (c)} \ldots )$. Then each element $x=(g,
\sigma)\in {\Gamma}_n$ gives rise to a partition-valued function
$( \rho (c))_{c \in \G_*} \in {\mathcal P} ( \G_*)$ such that
$\sum_{r, c} r m_r( \rho(c)) =n$. The partition-valued function $\rho =(
\rho(c))_{ c \in G_*} $ is called the {\em type} of $x$. It is
known that any two elements of ${\Gamma}_n$ are conjugate in
${\Gamma}_n$ if and only if they have the same type.

For $\lambda = (1^{m_1} 2^{m_2} \ldots )$, we define
$
  z_{\la } = \prod_{i \geq 1}i^{m_i}m_i!,
$ which is the order of the centralizer of an element of
cycle-type $\la $ in the symmetric
group $S_{|\la |}$. The order of the centralizer of
an element $x = (g, \sigma) \in {\Gamma}_n$ of the type $\rho=(
\rho(c))_{ c \in \G_*}$ is
\begin{eqnarray} \label{eq:centord}
Z_{\rho}=\prod_{c\in \G_*}z_{\rho(c)}\zeta_c^{l(\rho(c))}.
\end{eqnarray}

\begin{example}  \label{ex:running}
 Let $\G =\Z/2\Z =\{\pm\}$ with identity $+$,
and then $\Gn$ becomes the Weyl group of type
$B_n$ or $C_n$. Let $n=8$ and consider

\begin{eqnarray*}
 x &=& ( (+,+,-,+,+,+,-,-), (1,3)(2,4,7))\\
 y &=& ( (-,+,-,+,+,+,-,-), (1,7)(3,8)(4,5)).
\end{eqnarray*}
Then the cycle-product in $x$ associated to $(1,3)$ (resp. $(2,4,7)$)
is $-\cdot +=-$ (resp. $-\cdot +\cdot +=-$).
Thus the type of $x$ is given by:
the partition $(1,1)$ for the conjugacy class $+$
and the partition $(3,2,1)$ for the conjugacy class $-$.
Similarly, the type of $y$ is
given by: the partition $(2,2,2,1,1)$ for the conjugacy class $+$
and the empty partition $\emptyset$ for the conjugacy class $-$.
\end{example}
\subsection{The wreath products $\Gn$ and $\G_\infty$}

Let $\N$ be the set of positive integers.
The symmetric group $S_n$ acts on the set
$\{1,2 , \ldots, n\}$ by permutations and
$\Gn$ acts on $\{1,2 , \ldots, n\}$ by its projection onto $S_n$.
The wreath product $\Gn$ embeds
in $\G_{n+1}$ as the subgroup $\Gn \times 1$. The union
$\G_{\infty} = \cup_{n \ge 1} \Gn$
carries a natural group structure.
When $\G$ is trivial, $\G_\infty$ reduces to
$S_\infty  = \cup_{n \ge 1} S_n$.

For any set $Y$ of elements in $\G_{\infty}$, we define
a subset of $\N$
$$ \N(Y) =\{j \in \N \mid \sigma (j) \ne j
\text{  or } g_j \ne 1
\text{ for some }x=( (g_1, g_2, \ldots ),\sigma) \in Y \} .
$$
Clearly, we have $\N(Y) =\cup_{x\in Y}\N(x)$.
We denote by $| \N(Y)|$ the cardinality of $\N(Y)$.

\begin{example}
For $x, y$ in Example~\ref{ex:running}, we have
$\N(x) =\{1,2,3,4,7,8\}, \N(y) =\{1,3,4,5,7,8\}$
and $\N(x,y ) =\{1,2,3,4,5,7,8\}$.
\end{example}

If $(x_1, \ldots, x_r)$ and
$(y_1, \ldots, y_r)$ are two $r$-tuples of elements in $\G_\infty$,
we say that they are {\em conjugate} in $\G_\infty$ if
$y_i = z x_i z^{-1}$ $(i=1,\ldots, r)$ for some $z$
in $\G_\infty$. In this way, we divide $r$-tuples into
conjugate classes. For a given conjugate class $\mathcal C$ of
such $r$-tuples, we denote by $|\N(\mathcal C)|$ the cardinality of
$\N(x_1, \ldots, x_r)$  for any element
$(x_1, \ldots, x_r)$ in $\mathcal C$, as apparently
this number does not depends on the choice of such an element.
The next lemma is a generalization of
Lemma~2.1 in \cite{FH} for the symmetric group case.

\begin{lemma}   \label{lem:intersect}
 The intersection of a conjugate class $\mathcal C$ of $r$-tuples
of elements in $\G_\infty$ with $\overbrace{\Gn \times \cdots \times \Gn}^{r}$
is empty if $n < |\N(\mathcal C)|$, and is a
conjugate class of $\Gn$ if $n \ge |\N(\mathcal C)|$. The number
of $r$-tuples in this intersection is
$n (n-1)\ldots (n-|\N ({\mathcal C})|+1)/k(\mathcal C)$
where $k(\mathcal C) \in
\mathbb Q$ is a constant.
\end{lemma}

\begin{proof}
The first half of the lemma is obvious. Given an $r$-tuple $(x_1, \ldots, x_r)$
in this intersection, its stabilizer in $\Gn$ is
the direct product of its centralizer $H$
in the wreath product associated to the set $\N(x_1, \ldots, x_r)$
and the wreath
product associated to the complement of $\N(x_1, \ldots, x_r)$ in
$\{1,2 , \ldots, n\}$, so its order is
$A =k_1(\mathcal C) \cdot |\G|^{n -|\N ({\mathcal C})|}
 \cdot (n -|\N ({\mathcal C})|)!$,
where the constant $k_1(\mathcal C)$ is the order of $H$.
Therefore, the number of elements in the intersection (which is
a single conjugacy class) is $|\G|^n n!/A$, which can be recast
in the form as stated in the lemma, with
$k(\mathcal C) = k_1(\mathcal C)\cdot
|\G|^{- |\N ({\mathcal C})|}$.
\end{proof}
\subsection{The modified types}

Let $x$ be an element of $\Gn$ of type
$\rho =(\rho(c))_{c\in\G_*} \in \mathcal P_n(\G_*)$.
If we regard it as an element in $\G_{n+k}$
by the natural inclusion $\Gn \leq\G_{n+k}$, then $x$ has
the type $\rho\cup (1^k) \in \mathcal P_{n+k}(\G_*)$, where
$(\rho\cup (1^k)) (c) =\rho(c)$ for $c \neq c^0$ and
$(\rho\cup (1^k)) (c^0) =(\rho(c^0), 1, \ldots, 1) =\rho(c^0)\cup (1^k).$
It is convenient to introduce the {\em modified type}
of $x$ to be $\widetilde{\rho} \in
\mathcal P_{n- r}(\G_*)$, where $r =\ell(\rho(c^0))$, as follows:
$\widetilde{\rho}(c) =\rho(c)$ for $c \neq c^0$ and
$\widetilde{\rho}(c^0) = (\rho_1 -1, \ldots, \rho_r -1)$
if we write the partition ${\rho}(c^0) = (\rho_1, \ldots, \rho_r)$.
Two elements in $\G_\infty$ are conjugate if
and only if their modified types coincide.

Given $\mu \in \mathcal P(\G_*)$, we denote by
$\mathcal K_\mu$ the conjugacy class in $\G_\infty$ of elements
whose modified type is $\mu$.
If $x \in \mathcal K_\mu$, the {\em degree} $\|x\|$
of $x$ is defined to be the degree $\| \mu\|$ of its {\em modified}
type. We shall often denote by $(r)_c$,
where $r \in \Z_+$ and $c\in \G_*$,
an element in $\mathcal P(\G_*)$ consisting of a
single cycle of degree $r$ whose cycle-product lies in $c$.
We will say an element is of
{\em single-cycle modified type} $(r)_c$ if the modified type
of the element consists of a single cycle $(r)_c$.

\begin{example}
The modified type of $x$ in Example~\ref{ex:running}
is given by: the partition $\emptyset$
(resp. $(3,2,1)$) for the conjugacy class $+$ (resp. $-$).
The modified type of $y$ is given by:
the partition $(1,1,1)$ (resp. $\emptyset$) for $+$
(resp. $-$). Furthermore, $\|x\|=6$ and $\|y\| =3.$
\end{example}

\begin{remark}  \label{rem:basis}
 For each $n\ge 0$ and each $\mu \in \mathcal P(\G_*)$,
let $K_\mu (n)$ be the characteristic function of
the conjugacy class in $\Gn$ whose modified type
is $\mu$, i.e. the sum of all $\sigma \in \Gn \cap \mathcal K_\mu$.
Noting that
$| \N (\mathcal K_\mu)| =\Vert \mu \Vert +\ell(\mu (c^0))$,
we have by Lemma~\ref{lem:intersect}
that $K_\mu (n) \neq 0$ if and only if
$\Vert \mu \Vert +\ell(\mu (c^0)) \leq n$.
The nonzero  $K_\mu (n)$'s form a $\Z$-basis
for $R_\Z(\Gn)$.
\end{remark}

Given $g \in \G$ and $1 \le i \neq j$,
we denote by $(i \stackrel{g}{\rightarrow} j)$ the element
 $((g_1, g_2,\ldots), (i, j))$ in $\G_\infty$, where $(i,j) \in S_\infty$
is a transposition, $g_j=g, g_i =g^{-1}$, and
$g_k =1$ for $k \neq i,j$. Note that
$(i \stackrel{g}{\rightarrow} j) =(j \stackrel{g^{-1}}{\rightarrow} i)$
and it is of order $2.$
The elements in $\G_\infty$ of the form
$(i \stackrel{g}{\rightarrow} j)$, where
$g$ runs over $\G$ and $(i,j)$ (where $i<j$) runs over
all transpositions of $S_\infty$,
form the single conjugacy class whose cycle-products
are all $c^0$ and the partition corresponding to
$c^0$ is $(2,1,1,\ldots)$.
Clearly, any element $x$ in $\G_\infty$ can be written as a
product of elements
in $\mathcal K_{(1)_{c^0}}$ and elements of the form $h^{(i)} \in \G^{(i)}$,
$i \ge 1$. We call such a product a {\em reduced expression} for $x$
if $x$ cannot be written as a product of fewer such elements.
For a general element in $\G_\infty$, a reduced expression
can be constructed cycle-by-cycle.

\begin{example}  \label{ex:reduced}
Let $x = ((g_{i_1}, \ldots, g_{i_k}), (i_1 \ldots i_k))$
be of single-cycle modified type. The cycle-product
$p_x = g_k \cdots g_2 g_1$ lies in some
conjugacy class $c$.
The modified type of $x$ is $(k)_c$ if $p_x \neq 1$
and it is $(k-1)_c$ if $p_x = 1$.
A reduced expression for $x$ is given by

\begin{eqnarray*}
 x &=&
  \left\{
      \everymath{\displaystyle}
      \begin{array}{ll}
   (i_1 \stackrel{g_{i_2} }{\rightarrow} i_2)
   (i_2 \stackrel{g_{i_3}}{\rightarrow} i_3)
   \cdots ( i_{k-1} \stackrel{g_{i_k}}{\rightarrow} i_k) p_x^{(i_k)},
            &\text{if } p_x \neq 1,  \\
   (i_1 \stackrel{g_{i_2} }{\rightarrow} i_2)
   (i_2 \stackrel{g_{i_3}}{\rightarrow} i_3)
   \cdots ( i_{k-1} \stackrel{g_{i_k}}{\rightarrow} i_k),
                 & \text{if } p_x = 1.
      \end{array}
    \right.
\end{eqnarray*}
where $p_x^{(i_k)}$ denotes the $p_x$ in $\G^{(i_k)}$.
For later purpose, it is convenient to introduce
the following notation:
\begin{eqnarray*}
 \delta (p_x)  &=&
  \left\{
      \everymath{\displaystyle}
      \begin{array}{ll}
   1, &\text{if }  p_x \neq 1,  \\
   0, & \text{if } p_x = 1.
      \end{array}
    \right.
\end{eqnarray*}
Then, we have $\| x\| = (k-1) +\delta( p_x).$
\end{example}

\begin{lemma}   \label{lem:minimal}
Given an element $x \in \G_\infty$ of modified type $\la$,
the number of elements
appearing in a reduced expression of $x$ is $\Vert \la \Vert$.
\end{lemma}

\begin{proof}
Follows from the definitions of modified types and reduced expressions. (see the above example.)
\end{proof}

\begin{lemma}  \label{lem:ineq}
Let $x,y$ be elements in $\G_\infty$
of modified type $\la$ and $\mu$, and let
$x y$ be of modified type $\nu$. Then
$\Vert xy \Vert \le \Vert x \Vert + \Vert y \Vert$,
or equivalently,
$\Vert \nu \Vert \le \Vert \la \Vert + \Vert \mu \Vert.$
\end{lemma}

\begin{proof}
 Follows from Lemma~\ref{lem:minimal}.
\end{proof}
\subsection{Properties of the structure constants}

Given $\la, \mu \in  \mathcal P(\G_*)$, we write the convolution
product $K_ \la (n) \, K_\mu(n)$ in $R_\Z(\Gn)$ as a linear combination of
$K_\nu(n)$

$$K_\la (n) K_\mu(n) =\sum_\nu a_{\la\mu}^\nu (n) \; K_\nu(n)$$
where the structure constants $a_{\la\mu}^\nu (n) \in \Z_+$, and they are
zero unless $\Vert \nu \Vert \leq \Vert \la\Vert +\Vert \mu \Vert$ by
Lemma~\ref{lem:ineq}.
For $n \geq \Vert \nu \Vert + \ell( \nu(c^0))$,
the constant $a_{\la\mu}^\nu ( n)$
is uniquely determined (cf. Remark~\ref{rem:basis}).
For convenience, we sometimes denote by
$[A]B$ the coefficient of $A$ in $B$, and thus by
$[K_\nu(n)]K_\la (n) K_\mu(n)$ the structure
constant $a_{\la\mu}^\nu (n)$.

The following proposition generalizes Theorem~2.2 in \cite{FH}
in the symmetric group case.

\begin{proposition}  \label{prop:polydeg}
Let $\la,\mu,\nu \in \mathcal P(\G_*)$.
There is a unique polynomial $f_{\la\mu}^\nu (x)$
such that $f_{\la\mu}^\nu (n) =a_{\la\mu}^\nu (n)$
for all $n \ge \Vert \nu \Vert + \ell(\nu(c^0))$.
The degree of the  polynomial $f_{\la\mu}^\nu (x)$
is the greatest value of $|\N(\mathcal C)| - |\N(\mathcal K_\nu)|$
among all classes $\mathcal C$ of pairs $( x,y)$ such that
$x \in \mathcal K_\la, y \in \mathcal K_\mu,
xy \in \mathcal K_\nu.$
\end{proposition}

\begin{proof}
The set of pairs

\begin{eqnarray*}  \label{eq:cond}
P =\{(x,y) \in\G_\infty\times \G_\infty \mid
x \in \mathcal K_\la,  y \in \mathcal K_\mu,
xy \in \mathcal K_\nu \}
\end{eqnarray*}
is apparently invariant under the
simultaneous conjugation by $\G_\infty$,
so is a union of conjugate classes.

We claim that the set of
conjugate classes of such pairs in $P$ is finite. Indeed, if
$(x,y) \in P$, then
$|\N(x,y)| \leq |\N(x)| +|\N(y)|$
while $|\N(x)|$ and  $|\N(y)|$ are finite numbers
depending only on $\la$ and $\mu$
respectively. So by Lemma~\ref{lem:intersect},
each pair $(x,y) \in P$ is conjugate to some pair lying in
$\G_{|\N(x,y)|} \times \G_{|\N(x,y)|}$, which is
a finite set.

Let the conjugate classes of $P$ be
$\mathcal C_1, \mathcal C_2, \ldots, \mathcal C_r$.
Then, by Lemma~\ref{lem:intersect},
the number $T$ of pairs in the intersection
$P \cap (\Gn\times \Gn)$ is
$$T= \sum_{i=1}^r n(n-1) \cdots
(n- |\N(\mathcal C_i)| +1)/k(\mathcal C_i).$$

To obtain $a_{\la\mu}^\nu (n)$
for $n \ge \Vert \nu \Vert + \ell(\nu(c^0))$, we must
divide $T$ by the number of elements in
$\mathcal K_\nu \cap \Gn$, which is
$n(n-1) \cdots (n-|\N(\mathcal K_\nu)| +1)/k(\mathcal K_\nu),$
 again by Lemma~\ref{lem:intersect}.
Note that $|\N(\mathcal K_\nu)| = |\N(xy)| \leq
|\N (x,y)| =|\N(\mathcal C_i)|$ for $(x,y) \in \mathcal C_i$.
Thus,
$$a_{\la\mu}^\nu (n) =k(\mathcal K_\nu)\sum_{i=1}^r
(n-|\N(\mathcal K_\nu)|)(n-|\N(\mathcal K_\nu)| -1) \cdots
(n- |\N(\mathcal C_i)| +1)/k(\mathcal C_i)$$
is a polynomial of degree equal to the maximal value among
$|\N(\mathcal C_i)| -|\N(\mathcal K_\nu)|, 1\le i \le r$.
\end{proof}

\begin{remark}
 Proposition~5.1 of \cite{Wa4} (which is much more difficult
to establish) can be used to
derive the first part of Proposition~\ref{prop:polydeg}
and an upper bound of the degree
(but not the explicit formula for the degree).
\end{remark}

\begin{lemma} \label{lem:choice}
 Let $x = ( (g_1, g_2, \ldots), \sigma)$
be an element in $\G_\infty$, where $g_i \in \G^{(i)}$
and $ \sigma \in S_\infty$. If $\sigma (j)=i$ for $i\neq j$, then
we can find a reduced expression of $x$ which contains
the element $(j \stackrel{g_i}{\rightarrow} i)$.
\end{lemma}

\begin{proof}
One only needs to prove the lemma
in the case when $x$ consists of a single
nontrivial cycle, i.e., the cycle containing $i$
and $j$. We can write the cycle (by rotating cyclically if necessary)
so that the letter $i$ appears in the end
of the cycle.
In this case (with suitable change of letters),
a reduced expression as required has
been given in Example~\ref{ex:reduced}.
\end{proof}

\begin{proposition}  \label{prop:equaldeg}
Let $x,y \in \G_\infty$ be of modified types $\la$ and $\mu$, and
let $xy$ be of modified type $\nu$. If
$\Vert \nu \Vert = \Vert \la\Vert +\Vert \mu \Vert$,
then $\N(x,y) = \N(xy)$.
\end{proposition}

\begin{proof}
 Write $x = ( (g_1, g_2, \ldots), \sigma)$,
$y = ( (h_1, h_2, \ldots), \tau)$, and
$xy = ((f_1, f_2, \ldots), \sigma\tau)$,
where $g_i, h_i, f_i \in \G^{(i)}$.
Clearly, $\N(xy) \subset  \N(x,y)$.
If $i$ lies in $\N(x, y)$ but not in $\N(xy)$, then
by definition of the set $\N(\cdot)$, we have
$\sigma\tau (i) =i$ and $f_i =1$.
We have two possibilities:
a) $\tau (i) =j \neq i$, $\sigma (j)=i$,
and $g_i h_j =1$;
or b) $\sigma (i) =i$, $\tau (i) =i$,
and $g_i h_i =1$.

Let us first consider the case a).
By Lemma~\ref{lem:choice},
we can find a reduced expression of $x$ which contains
the element $(j \stackrel{g_i}{\rightarrow} i)$
and a reduced expression of $y$ which contains
the element $(i \stackrel{h_j}{\rightarrow} j)$.
Note that the elements
$(j \stackrel{g_i}{\rightarrow} i)$ and
$(i \stackrel{h_j}{\rightarrow} j)$ are
inverse to each other (and they are also equal)
thanks to $g_i h_j =1$.
Therefore, writing
$x y$ as the product of these two specified reduced
expressions, we can replace the part between (and which include)
$(j \stackrel{g_i}{\rightarrow} i)$ and
$(i \stackrel{h_j}{\rightarrow} j)$ by
another expression involving two less elements (by conjugation).
Thus, by Lemma~\ref{lem:minimal},
$\Vert \nu \Vert < \Vert \la\Vert +\Vert \mu \Vert$.

Let us now consider the case b). In fact, since $i$ is fixed
by both $\sigma$ and $\tau$, the element
$g_i$ (and thus $h_i$) is not equal to $1 \in \G$
by the assumption that $i$ lies in $\N(x, y)$.
As a reduced expression can be constructed cycle by cycle, we can
easily find a reduced expression for $x$ (resp. $y$) which satisfies
two requirements: 1) it contains
the element $g_i$ (resp. $h_i$) in $\G^{(i)}$;
2) $g_i$ (resp. $h_i$) commutes
with all the other elements appearing in this
reduced expression. Writing
$xy$ as the product of these two specified reduced
expressions, we can remove simultaneously $g_i$ and $h_i$
since $g_i h_i =1$.
Thus, again by Lemma~\ref{lem:minimal},
$\Vert \nu \Vert < \Vert \la\Vert +\Vert \mu \Vert$.
\end{proof}
\subsection{The Farahat-Higman ring}

Proposition~\ref{prop:polydeg}
and Proposition~\ref{prop:equaldeg} give us the following.

\begin{theorem}  \label{th:main}
Let $\la,\mu,\nu \in \mathcal P(\G_*)$. Then,
\begin{enumerate}
 \item
there is a unique polynomial $f_{\la\mu}^\nu (x)$
such that $f_{\la\mu}^\nu (n) =a_{\la\mu}^\nu (n)$
for all positive integers $n$.
 \item
the polynomial $f_{\la\mu}^\nu (x)$ is a constant if
$\Vert \nu \Vert = \Vert \la\Vert +\Vert \mu \Vert.$
\end{enumerate}
\end{theorem}

\begin{remark}
The number $a_{\la\mu}^\nu(n)$
can take any value for $n <\Vert \nu \Vert+ \ell( \nu(c^0))$
since $K_\nu(n) =0$ by Remark~\ref{rem:basis}.
In the above theorem, it is understood that
$a_{\la\mu}^\nu (n)$ for
$n <\Vert \nu \Vert+ \ell( \nu(c^0))$
is chosen appropriately
(i.e. equal to the value of the polynomial $f_{\la\mu}^\nu (x)$ at
$x=n$).
\end{remark}

\begin{remark}
The inverse of Theorem~\ref{th:main}~(2) is not true even
in the symmetric group case (i.e. when $\G$ is trivial),
in contrast to the claim in Ex.~24, pp. 131, \cite{Mac}. For
example, consider the conjugacy classes
$\mathcal K_\la, \mathcal K_\mu$ and $\mathcal K_\nu$
of modified types
$\la =(1), \mu =(3)$ and $\nu= (1,1)$.
Clearly, $\Vert \nu \Vert < \Vert \la\Vert +\Vert \mu \Vert.$
For any $\sigma \in \mathcal K_\la,\tau\in \mathcal K_\mu$ such
that $\sigma\tau \in \mathcal K_\nu$, it is easy
to see that $\sigma$, $\tau$ and $\sigma\tau $
have to be of the form $\sigma =(i,k)$, $\tau =(i,j,k,l)$
and $\sigma\tau =(i,j)(k,l)$, where
$i,j,k,l$ are distinct positive integers. If follows that
$\N(\sigma\tau) = \N(\sigma,\tau) =\{i,j,k,l\}$.
Thus, by Proposition~\ref{prop:polydeg}, the polynomial
$f_{\la\mu}^\nu $ is a constant. One can further show
that $f_{\la\mu}^\nu $ is actually $2$ in this case,
by showing that there is a single conjugate classes
of such pairs $(\sigma,\tau)$ and then using the formula for
$a_{\la\mu}^\nu (n)$ in the proof of
Proposition~\ref{prop:polydeg}.
\end{remark}

Denote by $R$ the subring of the polynomial ring
$\mathbb Q [t]$ consisting of polynomials that
take integer values at all integers. Thanks to
Theorem~\ref{th:main}, we can introduce
a commutative $R$-ring $\mathcal F_\G$ with $R$-basis
$(K_\la)$ indexed by $\la \in \mathcal P(\G_*)$
and multiplication given by
$$K_\la  K_\mu =\sum_\nu f_{\la\mu}^\nu  K_\nu$$
where the sum is over $\nu \in \mathcal P(\G_*)$
such that $\Vert \nu \Vert \leq \Vert \la\Vert +\Vert \mu \Vert,$
and $f_{\la\mu}^\nu \in R$ takes the value
$a_{\la\mu}^\nu (n)$ at a positive integer $n$.

The $R$-ring $\mathcal F_\G$ determines
the class algebra $R_\Z(\Gn)$ for each $n$ by
the natural surjective ring homomorphism
$\mathcal F_\G \longrightarrow R_\Z(\Gn)$ sending
$K_\la$ to $K_\la (n)$ and $f_{\la\mu}^\nu$
to $a_{\la\mu}^\nu (n)$.

The assignment of degree $\Vert \la\Vert$ to $K_\la$
provides a filtration on the ring $\mathcal F_\G$,
thanks to Lemma~\ref{lem:ineq}. The associated graded
ring, denoted by $\mathcal G_{R,\G}$, has the multiplication

\begin{eqnarray} \label{eq:const}
K_\la \,  K_\mu
 =\sum_{\Vert \nu \Vert = \Vert \la\Vert +\Vert \mu \Vert}
 a_{\la\mu}^\nu  K_\nu .
\end{eqnarray}
By abuse of notations, we have kept using the same symbols
for elements in the ring $\mathcal F_\G$ and the associated graded ring.
By Theorem~\ref{th:main}, the structure constants for
$\mathcal G_{R,\G}$ are {\em non-negative integers}.
Thus, we have reached the following.

\begin{theorem}
If $\mathcal G_{\G}$ denotes
the ring over $\Z$ with basis $(K_\la)$ whose
multiplication is given by (\ref{eq:const}),
then the graded $R$-ring $\mathcal G_{R,\G}$
associated to $\mathcal F_\G$ is
isomorphic to $ R \otimes_\Z \mathcal G_{\G}$.
\end{theorem}
The ring $\mathcal G_\G$
will be called the {\em Farahat-Higman ring} (or {\em FH-ring}) associated
to $\G$, as our results specialize in the symmetric group
case to the ones  obtained by Farahat-Higman \cite{FH}.

We can also take the filtration for each $n$ directly.
Recall that those $K_\mu (n)$ indexed by $\mu \in \mathcal P(\G_*)$
satisfying
$\Vert \mu \Vert +\ell(\mu (c^0)) \leq n$ form a basis
for $\Rgn$. The assignment of degree
$\Vert \mu\Vert$ to $K_\mu (n)$ provides $R_\Z(\Gn)$
a filtered ring structure, whose
associated graded ring will be denoted by $\mathcal G_\G (n)$.
The next theorem is essentially a reformulation
of Theorem~\ref{th:main} for $\Gn$.

\begin{theorem}  \label{th:constant}
The graded $\Z$-ring $\mathcal G_\G (n)$ has a multiplication given by

$$K_\la (n) \, K_\mu(n)
=\sum_{\Vert \nu \Vert = \Vert \la\Vert +\Vert \mu \Vert}
 a_{\la\mu}^\nu  K_\nu(n)$$
where the structure constants $a_{\la\mu}^\nu $
are nonnegative integers independent of $n$.
\end{theorem}
Clearly, the ring $\mathcal G_\G$ determines
the family of rings $\{ \mathcal G_\G (n) \}_{n\ge 1}$,
and vice versa. There is a natural surjective ring
homomorphism $\mathcal G_\G \rightarrow \mathcal G_\G (n)$
for each $n$ which is compatible with the surjective ring
homomorphism $\mathcal G_\G (n+1) \rightarrow \mathcal G_\G (n)$
by restriction.
When $\G$ is the trivial group and thus
$\Gn$ reduces to the symmetric group $S_n$,
we omit the subscript $\G$ from the notations for the rings
$\mathcal F_\G$, $\mathcal G_\G$, $\mathcal G_\G(n)$, etc.

\begin{remark}
 In the symmetric group case (i.e. when $\G$ is trivial),
the above results are due to Farahat-Higman \cite{FH}.
We have adopted however the elegant presentation of Macdonald
(cf. Ex.~24, \cite{Mac}, pp. 131). Macdonald in addition found
a symmetric function interpretation for the FH-ring
$\mathcal G =\mathcal G_1$ (for $\G=1$).
We expect that our FH-ring $\mathcal G_\G$ also affords a natural
symmetric function interpretation.
\end{remark}
\section{Computations in the Farahat-Higman ring}
\label{sec:computation}

In this section, we will compute the formulas
for the multiplication by a conjugacy class
of single-cycle modified type
in the FH-ring $\mathcal G_\G$, and then
derive various consequences.

\subsection{Constraints from the filtration}

Let $x, y \in \G_\infty$, where we assume that
$y$ is of single-cycle modified type, i.e., $y$ is of the form
$y =((g_{i_1}, \ldots, g_{i_r}), (i_1,  \ldots, i_r))$.
We will say that a cycle
of $x$ is {\em relevant} to $y$ if
it contains at least one of $i_1, \ldots, i_r$.

\begin{lemma}  \label{lem:sameprod}
Let $y =((h_{i_1}, \ldots, h_{i_r}), (i_1, \ldots, i_r))$ be of
single-cycle modified type $(r >2)$. Then, we can decompose $y$ as
a product $y_1 y_2$, where  $y_1 =(f, (i_1, \ldots, i_{r-1}))$ and
$y_2 =((a^{-1},a), (i_{r-1}, i_r)) \in \mathcal K_{(1)_{c^0}}$,
for some $a \in\G, f \in \G^{r-1}$, such that $\|y\| =\| y_1\|
+\|y_2\|$. For any such decomposition, the cycle-products of $y$
and $y_1$ are conjugate in $\G$.
\end{lemma}

\begin{proof}
Fix  $y_2$ by choosing an $a \in\G$. Then
one easily checks that $y_1 =yy_2^{-1}$
can be recast in the form in the lemma.
The second statement
follows from a straightforward computation.
\end{proof}

The next key proposition is a purely wreath product
phenomenon, which reduces to almost triviality
in the symmetric group case. It
indicates that the multiplication
in the FH-ring is rigidly constrained by the filtration
given by $\| \cdot \|$. Below  $\| \cdot \|$ also
applies to elements in $S_\infty$ by regarding
$S_\infty$ as the natural subgroup of $\G_\infty$.

\begin{proposition} \label{prop:constraint}
 Let $y =((h_{i_1}, \ldots, h_{i_r}), \tau_r)$ be of
single-cycle modified type, where
$\tau_r =(i_1 i_2 \ldots i_r)$. Let $x =(g, \sigma) \in \G_\infty$
be such that $\|x\| + \|y\| = \|xy \|$. Then
one of the following possibilities occurs:

\begin{enumerate}
 \item $\|\sigma \tau_r\|= \| \sigma\| +\|\tau_r \|$.
In this case, there are $r$ cycles of $x$ relevant to $y$,
i.e., each $i_s$ $(s=1,\ldots, r)$ belongs
to a distinct cycle in $x$;
among the cycle of $y$ and all $r$ cycles of $x$ relevant to $y$,
at most one cycle-product is not $1$. When
all such cycle-products (resp. all but one, say $p$)
are $1$, then the relevant
cycle-product in $xy$ is also $1$ (resp. belong to the
conjugacy class of $p$).

\item $\|\sigma \tau_r\|= \| \sigma\| +\|\tau_r \| -2 $.
In this case,
there are $(r-1)$ cycles of $x$ relevant to $y$;
all the cycle-products of the relevant cycles in $x$ are $1$
and so is
the cycle-product of $y$. Assume that $i_s$ and $ i_t$
are contained in the same relevant cycle of $x$,
then $i_s$ and $i_t$ are contained in distinct cycles
of $xy$ whose cycle-products lie in opposite
conjugacy classes.
\end{enumerate}
\end{proposition}

\begin{proof}
Let us denote $n =| \N (x,y)|$, and regard
$x$ and $y$ as elements in $\Gn$.
Since $\| \cdot \|$ is additive for disjoint cycles,
we may assume that all
the cycles of $x$ are relevant to $y$
without loss of generality. If we denote by $k$ the
number of (relevant) cycles appearing in $x$, then $k \le r$.

Denote by $\Delta_x$ the number of cycles in $x$
whose cycle-products are not $1$, and by $p_y$
the cycle-product of $y$ as usual.
It follows from the definitions that
$\|x\| = (n-k) +\Delta_x$, and
$\|y\| =(r-1) + \delta (p_y)$. Thus, we have

\begin{eqnarray}
\|x\| +\|y\| &=& n +(r-k)-1 + \Delta_x + \delta (p_y)
  \label{eq:ident1}\\
\|x y\| &\le& n
  \label{eq:ineq2}\\
\|x\| +\|y\| &=& \|x y\|
  \label{eq:assump}
\end{eqnarray}
where (\ref{eq:ineq2}) holds since $xy \in\Gn$,
and (\ref{eq:assump}) by the assumption.

If follows from the consistency of
(\ref{eq:ident1}, \ref{eq:ineq2}, \ref{eq:assump})
and $k \le r$
that either (1) $k=r$ or (2) $k=r-1$.

 Let us first consider the case (1), i.e. $k=r$.
There are $r$ cycles of $x$ relevant to $y$,
i.e., each $i_s$ $(s=1,\ldots, r)$ belongs
to a distinct cycle in $x$. It follows that
$xy$ consists of a single $n$-cycle
(see e.g. the formula (\ref{eq:link}) below).
In particular, $\| \sigma \tau_r \| =n-1 =\| \sigma \| +\|\tau_r \|$.
By the comparison of
(\ref{eq:ident1}, \ref{eq:ineq2}, \ref{eq:assump}),
we have two subcases:
(1a)  $\|x y\|=n -1$, and  $\Delta_x =\delta (p_y)= 0$;
(1b) $\|x y\|=n$, and then one of
$\Delta_x$ and $\delta (p_y)$ is 0 while the other is $1$.

In the subcase (1a),
the cycle of $y$ and the $r$ cycles of $x$
all have cycle-product equal to $1$.
In the subcase (1b),
among the cycle of $y$ and the $r$ cycles of $x$,
all cycle-products but one
are equal to 1. Let us denote the non-identity cycle-product
by $p$.

It remains to check the last statement in part (1) of the
proposition. For $r=1$ or $2$, it is straightforward
to check by hand.
For $r >2$, we can write $y$ as a product $y_1 y_2$,
where $y_1 =(h, (i_1, \ldots, i_{r-1}))$
and $y_2 =((a^{-1},a), (i_{r-1}, i_r))$,
by Lemma~\ref{lem:sameprod}. Now, an obvious
induction proves the last statement of the part (1) of the proposition.

Let us now consider the case (2),
i.e. $k =r-1$. It immediately follows
that there are $(r-1)$ cycles of $x$ relevant to $y$.
Furthermore, $\sigma \tau_r$ consists of
exactly two cycles (see e.g. the formula (\ref{eq:mult}) below). Thus,
$\|\sigma \tau_r\| =n-2$ and this coincides with
$\|\sigma\| + \|\tau_r\| -2 = (n -k)+(r-1) -2=n-2.$

It follows again from the consistency of
(\ref{eq:ident1}, \ref{eq:ineq2}, \ref{eq:assump})
that $\Delta_x = \delta (p_y) =0$
(and $\|xy \|=n$), i.e., the cycle of $y$
and the $r$ cycles of $x$ all
have cycle-product $1$.

It remains to check the last statement in part (2) of the
proposition. For $r=2$, it is straightforward
to check by hand.

Now let $r>2$. There are exactly two $i$'s
among $i_1, \ldots, i_r$ which
lie in an identical cycle of $x$. We
can assume one of them is $i_r$ by rotating cyclically
$(i_1, \ldots, i_r)$ if necessary.
Then, by Lemma~\ref{lem:sameprod},
we can write $y =y_1 y_2$, where
$y_1 = (f, (i_1, \ldots, i_{r-1}))$ and
$y_2 = ((a^{-1}, a), (i_{r-1}, i_r))$
such that $\|y\| =\|y_1\|+\|y_2\|$.
Together with the identity $\|x\| + \|y\| = \|x y\|$
and the sub-multiplicative property of $\| \cdot \|$
(cf. Lemma~\ref{lem:ineq}),
we have $\|x\| + \|y_1\| = \|x y_1\|$
and $\|xy_1\| + \|y_2\| = \|x y\|$.

Since $i_1, \ldots, i_r$ lie in the
$(r-1)$ distinct cycles of $x$,
the product $xy_1$ is of single-cycle modified type.
Thus, we are in the setup of (1) when considering
the product $x$ and $y_1$; and in the setup
of (2) when considering the product $(xy_1)$ and $y_2$.
By applying (2) for $r=2$
(which we have checked) to the product $(xy_1)$ and $y_2$,
we conclude that the cycle-product of
$(xy_1)$ is $1$.
Now, by applying (1) to the product $x$ and $y_1$,
we finish the proof of
the part (2) of the proposition for general $r$.
\end{proof}
\subsection{Multiplication with a conjugacy class of
single cycles}

We are ready to compute the explicit formulas for
multiplications in the FH-ring $\mathcal G_\G$
with a conjugacy class of single-cycle
modified type. Note that $K_{(0)_{c^0}} =K_\emptyset$ is the identity of
the FH-ring $\mathcal G_\G$.

\begin{theorem}  \label{th:singlecyc1}
 Let $ \la\in \mathcal P(\G_*)$, $r \geq 0$,
 and $c\in \G_* \backslash c^0$. Then,
 the product $K_\la K_{(r)_{c}}$ in the FH-ring
$\mathcal G_\G$ is

 \begin{eqnarray*}
 K_\la K_{(r)_{c}} =
 \sum_{\mu}
    \frac{\left( m_{\|\mu \| +r}( \la (c)) +1 \right)
 \cdot (\|\mu \| +r) \cdot (r-1)!}{(r -\ell (\mu))! \prod_{i\ge 1}m_i(\mu)!}
 K_{\la \cup (\|\mu \| +r)_{c} -\mu}
 \end{eqnarray*}
 summed over
 $\mu=\mu(c^0)\subset \la(c^0)$ such that $\ell (\mu) \le r$,
where we denote $\mu = (1^{m_1(\mu)} 2^{m_2(\mu)} \ldots)$.
\end{theorem}

\begin{proof}
 Let $y =((h_{i_1}, \ldots, h_{i_r}), (i_1 \ldots i_r))$
be of type $(r)_c$, i.e. $h_{i_r} \cdots h_{i_1} \in c$.
Let $x =(g, \sigma)$ be in the class $\mathcal K_\la$
such that $\|xy\| = \|x\|+\|y\|.$
By Proposition~\ref{prop:constraint}~(1), the
numbers $i_1, \ldots, i_r$
belong to $r$ distinct cycles of $\sigma$, which we denote
by $(\ldots j_a i_a)$, $a=1, \ldots, r$; the relevant
cycle-products for $x$ are all 1.

Note that
 \begin{eqnarray}  \label{eq:link}
 (\ldots j_1 i_1) (\ldots j_2 i_2) \cdots (\ldots j_r i_r)
\cdot (i_1, \ldots i_r)
= (\ldots j_1 i_1 \ldots j_2 i_2 \ldots j_r i_r).
\end{eqnarray}
Assume that $m_0(\mu)$ of the cycles of
$x$ relevant to $y$
(i.e. those which we have used on the left-hand side
of (\ref{eq:link})) are 1-cycles,
$m_1(\mu)$ 2-cycles,
$m_2 (\mu)$ 3-cycles, etc., of $x$. Then,
$\mu =\mu (c^0) =(1^{m_1(\mu)}2^{m_2 (\mu)} \ldots)$
is the modified type of the product of cycles in $x$ relevant to $y$,
and $m_0(\mu) =r- \ell (\mu)$. By Proposition~\ref{prop:constraint}~(1),
the cycle-product associated to the cycle
$ (\ldots j_1 i_1 \ldots j_2 i_2 \ldots j_r i_r)$ of $x$
lies in $c$. Then $xy$ belongs to
$\mathcal K_{\la \cup (\| \mu\| +r)_c -\mu} $.
This implies that the conjugacy classes on the right-hand
side of the formula in the theorem
are the right ones.

Next, we determine the coefficient of
$K_{\la \cup (\| \mu\| +r)_c -\mu} $.

Fix $z \in \mathcal K_{\la \cup (\| \mu\| +r)_c -\mu}$ and
$x_1 \in  \mathcal K_\mu$. This coefficient is
the number of ways a cycle of modified type $(\| \mu\| +r)_c$
can be chosen from $z$, and bracketted
to form a product of  $m_0(\mu)$ cycles of
modified type $(0)_{c^0}$, $m_1(\mu)$ cycles of
modified type $(1)_{c^0}$, $m_2(\mu)$ cycles of
modified type $(2)_{c^0}$, etc., and a cycle of
modified type $(r)_{c}$. The cycle of modified type $(\| \mu\| +r)_c$
can be chosen in $m_{\| \mu\| +r} (\la (c)) +1$ ways.
The number of brackettings, which no longer depends
on $\la$ but only on $\mu$ and $(r)_c$, is the
coefficient  $[K_{(\| \mu\| +r)_c}] K_{\mu} K_r (c)$.
By Theorem~\ref{th:main}, this coefficient is
a constant, so we can calculate it as the coefficient
$[K_{(n)_c}(n)] K_{\mu}(n) K_{(r)_c}(n)$, where
$n =\| \mu\| +r$.

The number of elements $u$ in
$K_{\mu}(n)$ is, by (\ref{eq:centord}), given by

\begin{eqnarray} \label{eq:conjord}
 \frac{|\G|^n n!}{ |\G|^r 2^{m_1(\mu)} 3^{m_2(\mu)} \cdots
 (r -\ell(\mu))! \prod_{i \ge 1} m_i(\mu)!}.
\end{eqnarray}

To find a cycle $y$ of modified type $(r)_{c}$
such that $u y$ is an $n$-cycle (which will automatically be
of modified type $(n)_c$ by Lemma~\ref{lem:sameprod}),
we need to choose a number from each of the $r$ cycles
of $u$ as well as $r$ elements from $\G$,
and then arrange them into a cycle of modified type $(r)_{c}$.
The number of choices here is

\begin{eqnarray} \label{eq:match}
 2^{m_1(\mu)} 3^{m_2(\mu)} \cdots  (r-1)! |\G|^r/\zeta_c.
\end{eqnarray}
Finally, to obtain the required coefficient,
we need to divide the product of (\ref{eq:conjord}) and
(\ref{eq:match}) by the number of elements in
$K_{(n)_c}(n)$ which is $|\G|^n n!/ n \zeta_c$.
Remembering $n =\| \mu\| +r$, we have proved the theorem.
\end{proof}

Given a partition $\mu =\mu (c^0) =(\mu_1, \mu_2, \ldots, \mu_k)$
where $\mu_1 \geq \mu_2 \geq \ldots \geq \mu_k \geq 1$
and $k =\ell(\mu)$. Assume $r \geq \ell(\mu)$.
we denote by $\tilde{\mu}$ the partition
$\tilde{\mu}=(\mu_1+1, \mu_2+1, \ldots, \mu_k+1, 1, \ldots, 1)$
of $|\mu| +r$ which consists of $r$ parts.
There are $r!$ compositions (counted with multiplicities)
obtained from rearranging the $r$ parts of $\tilde{\mu}$.
Given $s>0$, $\mu$ and $r \geq \ell(\mu)$,
we denote by $q(\mu,r,s)$ (resp. $p(\mu,r,s)$)
the number of all compositions
$(\tilde{\mu}_{k_1},\tilde{\mu}_{k_2}, \ldots, \tilde{\mu}_{k_r})$
associated to $\tilde{\mu}$,
counted with (resp.  without) multiplicities,
such that $s$ is not one of the $r$ numbers
$\tilde{\mu}_{k_1} + \tilde{\mu}_{k_2}
 + \ldots +\tilde{\mu}_{k_a} (a=1, \ldots, r).$
If we write $\mu =\mu (c^0)
=(1^{m_1(\mu)} 2^{m_2(\mu)}\ldots)$, then
$\tilde{\mu} =(1^{r -\ell(\mu)} 2^{m_1(\mu)} 3^{m_2(\mu)}\ldots)$,
and

$$q (\mu,r,s) =
{(r-\ell (\mu))! \prod_{i\ge 1}m_i(\mu)!} \cdot p(\mu,r,s).$$
Clearly, $p(\mu,r,s)$ is zero unless $0< s < |\mu| +r$.
 We mention that there is a symmetry:
$p(\mu,r, s) =p(\mu,r, |\mu| +r -s),$ and
$\sum_{s} q (\mu,r, s) =\|\mu \| \cdot r!.$ But we
do not need these properties in this paper.

\begin{theorem}  \label{th:singlecyc2}
Let $ \la\in \mathcal P(\G_*)$ and $r \geq 0$. Then,
the product $K_\la K_{(r)_{c^0}}$ in the FH-ring
$\mathcal G_\G$ is
  \begin{eqnarray*}
 K_\la K_{(r)_{c^0}}
 &=&
 \sum_{I_1}
\frac{\left(m_{\|\mu \| +r} (\la (c^0)) +1 \right)
 \cdot (\|\mu \| +r +1) \cdot r!}{ (r+1-\ell (\mu))!
 \prod_{i\ge 1} m_i(\mu)!} K_{\la \cup (\|\mu \| +r)_{c^0} -\mu}
   \\
 & +&
   \sum_{I_2}
 \frac{\left(m_{\|\mu \| +r+k} (\la (c)) +1 \right)
 \cdot (\|\mu \| +r+k) \cdot r!}{(r-\ell (\mu))! \prod_{i\ge 1}m_i(\mu)!}
 K_{\la \cup (\|\mu \| +r +k)_c -(\mu \cup (k)_c)}
   \\
 & +&
  \sum_{I_3}
    s_1s_2 \zeta_c p (\mu,r,s_1)
     (m_{s_1}(\la(c)) +1)(m_{s_2}(\la(\bar{c})) +1)
    K_{\la \cup (s_1)_c \cup (s_2)_{\bar{c}}-\mu}\\
 & +&
  \sum_{I_4}
     s^2 \zeta_c p(\mu,r,s)
     (m_{s}(\la(c)) +1) (m_{s} (\la(\bar{c})) +1)
    K_{\la \cup (s)_c \cup (s)_{\bar{c}}-\mu}\\
 & +&
  \sum_{I_5}
   s^2 \zeta_c p(\mu,r,s)
    { m_{s}(\la(c)) +2 \choose 2}
    K_{\la  \cup  (s)_c \cup  (s)_c-\mu}
 \end{eqnarray*}
 where  $\mu = (1^{m_1(\mu)} 2^{m_2(\mu)} \ldots)$,
$I_1 =\{\mu \mid \mu =\mu(c^0) \subset \la(c^0) ,
           \ell (\mu) \le r +1 \}$,
$I_2 =\{ ( c,\mu, k) \mid
 c \neq c^0, \mu =\mu(c^0) \subset \la(c^0),
 \ell (\mu) \le r,  (k)_c \subset \la(c), k\geq 1 \}$,
$I_3 =\{( c,\mu, s_1,s_2) \mid c \neq c^0 ,
   \mu =\mu(c^0) \subset \la(c^0),
           \ell (\mu) \leq r ,
  0 <s_1 <s_2, s_1 + s_2=\|\mu \| +r \}$,
$I_4 =\{( \{c, \bar{c} \},\mu, s) \mid
  c \neq \bar{c},
  \mu =\mu(c^0) \subset \la(c^0),  \ell (\mu) \le r,\,
  2s =\|\mu \| +r \}$, and
$I_5 =\{( c,\mu, s) \mid
  c = \bar{c}, c \neq c^0,
  \mu =\mu(c^0) \subset \la(c^0),  \ell (\mu) \le r,\,
  2s =\|\mu \| +r \}$.
 Here $\{c, \bar{c} \}$ denotes an unordered pair
of conjugacy classes in $\G_*$.

\end{theorem}

\begin{proof}
 Let $y =((h_{i_1}, \ldots, h_{i_r}, h_i), (i_1 \ldots i_r i))$
be of type $(r)_{c^0}$, i.e. $h_i h_{i_r} \cdots h_{i_1} =1$.
Let $x =(g, \sigma)$ be in the class $\mathcal K_\la$
such that $\|x \| +\|y \| =\|x y\|.$
Both cases of Proposition~\ref{prop:constraint} occur
when multiplying $x$ and $y$. We will analyze them one by one.
We keep in mind the strategy in the proof of
Theorem~\ref{th:singlecyc1}.

(1). The numbers $i_1, \ldots, i_r, i$
belong to $(r+1)$ distinct cycles of $\sigma$, which we denote
by $(\ldots j i)$ and
$(\ldots j_a i_a)$, $a=1, \ldots, r$. We divide this into
two subcases according to Proposition~\ref{prop:constraint}~(1):

(1a) all the relevant
cycle-products for $x$ are 1; (1b)  all
the relevant cycle-products but one (which
lies in $c \in \G_* \backslash c^0$)
for $x$ are 1.

In the case (1a), we have the same computation as
in the proof of Theorem~\ref{th:singlecyc1}, which gives us
the first summand on the right-hand side of the formula
in the theorem. The main difference we should keep in
mind is that $y$ is an $(r+1)$-cycle (not an $r$-cycle)
since $y$ is of {\em modified} type $(r)_{c^0}$.
We omit the detail here.

In the case (1b),
we use from $x$ the relevant cycles of modified
type $\mu =\mu (c^0)$ (which is contained in $\la (c^0)$)
as well as a cycle of modified type $(k)_c$
(which is contained in $\la (c)$). When multiplying
with $y$, the same strategy as in the proof of
Theorem~\ref{th:singlecyc1} applies, as we sketch below.

The number of brackettings is the
coefficient
$[K_{(n)_c}(n)] K_{\mu\cup (k)_c}(n) K_{(r)_{c^0}}(n)$, where
$n =\| \mu\| +r +k$.

The number of elements $u$ in $K_{\mu\cup (k)_c}(n)$ is

\begin{eqnarray} \label{eq:conum}
 \frac{|\G|^n n! }{|\G|^r \zeta_c k\cdot 2^{m_1(\mu)} 3^{m_2(\mu)} \cdots
 (r - \ell(\mu))! \prod_{i\geq 1}m_i(\mu)!}.
\end{eqnarray}

To find a cycle $y$ of modified type $(r)_{c^0}$
such that $u y$ is an $n$-cycle (which will automatically be
of modified type $(n)_c$ by Lemma~\ref{lem:sameprod}),
we need to choose a number from each of the $r$ relevant
cycles of $u$ (whose cycle-product is 1),
a number from the cycle of modified type $(k)_c$,
as well as $(r+1)$ elements from $\G$,
and then arrange into a cycle of modified type $(r)_{c^0}$.
The number of choices of a cycle $y$ of modified type $(r)_{c^0}$
such that $u y$ is an $n$-cycle is

\begin{eqnarray} \label{eq:numb}
 k \cdot 2^{m_1(\mu)} 3^{m_2(\mu)} \cdots r! |\G|^r.
\end{eqnarray}

The product of (\ref{eq:conum}) and (\ref{eq:numb})
divided by the number of elements in
$K_{(n)_c}(n)$ gives us the coefficient
$[K_{(n)_c}(n)] K_{\mu\cup (k)_c}(n) K_{(r)_{c^0}}(n)$.
Noting  the number of elements in
$K_{(n)_c}(n)$ is $|\G|^n n!/ n \zeta_c$
and remembering $n =\| \mu\| +r +k$, we obtain
the coefficient in the
second term on the right-hand side of the formula
in the theorem.

(2). Assume that $i_t (t=1, \ldots, r)$ belong
to distinct cycles, say
$(\ldots j_t i_t)$, of $x$, and we also assume that
$i_a$ ($a \leq r$) and $i_r$ lie in the same cycle of $x$.
By Proposition~\ref{prop:constraint} (2),
the cycle-product of $y$ and all the relevant
cycle-products for $x$ are 1.

Note that
\begin{eqnarray}
 && (\ldots j_1 i_1)  \cdots (\ldots j_{a-1} i_{a-1})
 (\ldots j i \ldots j_a i_a)
 \cdots  (\ldots j_{r} i_{r})
 \cdot (i_1 \ldots i_r i) \nonumber \\
 &=& (\ldots j_1 i_1 \ldots j_{a-1} i_{a-1}
  \ldots j i \ldots j_a i_a\ldots j_{r} i_{r})
  (i_{r} i)  \nonumber \\
 &=& (\ldots j_1 i_1 \ldots j_{a-1} i_{a-1} \ldots j i)
  (\ldots j_a i_a\ldots j_{r} i_{r})
  \label{eq:mult}
\end{eqnarray}
The cycle-products in $xy$ of the
two cycles, say $s_1$-cycle and $s_2$-cycle,
on the right-hand side of (\ref{eq:mult})
belong to opposite conjugacy classes,
by Proposition~\ref{prop:constraint} (2).
We denote the modified type of
the right-hand side of (\ref{eq:mult})
to be $(s_1)_c \cup (s_2)_{\bar{c}}$.
Let us assume that $m_0(\mu)$ of the cycles of
$x$ (relevant to $y$) are 1-cycles, $m_1(\mu)$ 2-cycles,
$m_2 (\mu)$ 3-cycles, etc. We denote by
$\mu =\mu (c^0) =(1^{m_1(\mu)}2^{m_2 (\mu)} \ldots)$
the modified type of the product of cycles in $x$ relevant to $y$,
and hence $m_0(\mu) =r- \ell (\mu)$.

Thus, $xy$ belongs to the conjugacy class
$\mathcal K_{\la \cup (s_1)_c \cup (s_2)_{\bar{c}} -\mu }$
which give rise to
the third, fourth, and fifth terms of the right-hand
side of the formula in the theorem.
So we see that the conjugacy classes
appearing on the right-hand side of the formula
in the theorem are indeed the right ones.

It remains to compute the coefficient
$N := [K_{\la \cup (s_1)_c \cup (s_2)_{\bar{c}} -\mu }]
K_{\la} K_{(r)_{c^0}}$,
where $s_1 +s_2 =\| \mu\| +r$.

We denote by $B$ the coefficient
$[K_{(s_1)_c \cup (s_2)_{\bar{c}} }] K_{\mu} K_{(r)_{c^0}}$.
Denote by $A_1$ the number of ways of
cycles of modified type $(s_1)_c \cup (s_2)_{\bar{c}}$
can be chosen from a fixed element in
$\mathcal K_{\la \cup (s_1)_c\cup (s_2)_{\bar{c}} -\mu }$.
Clearly, we have $N =A_1 \cdot B$.
Note that $B$ does not depend on $\lambda$
and only depends on $\mu,r$ and $s_1$.
By Theorem~\ref{th:main},
the number $B$ can be computed as the coefficient
$[K_{(s_1)_c \cup (s_2)_{\bar{c}} }(n)] K_{\mu}(n) K_{(r)_{c^0}}(n)$,
where $n =\| \mu\| +r.$

We first easily have
 \begin{eqnarray*}
 A_1
  &=&
  \left\{
      \everymath{\displaystyle}
      \begin{array}{ll}
    {m_{s}(\la(c)) +2 \choose 2}
  , &\text{if }  c=\bar{c} \text{ and } s_1=s_2=s,  \\
   (m_{s_1}(\la(c)) +1)(m_{s_2}(\la(\bar{c})) +1),
   &  \text{otherwise}.
      \end{array}
    \right.
\end{eqnarray*}

The number $A_2$ of elements $u$ in $K_\mu(n)$ is,
by (\ref{eq:centord}), given by

\begin{eqnarray*}
A_2 =
\frac{|\G|^n n! }{|\G|^r 2^{m_1(\mu)} 3^{m_2(\mu)} \cdots
 (r -\ell (\mu))! \prod_{i \geq 1} m_i(\mu)!}.
\end{eqnarray*}

To find a cycle $y$ of modified type $(r)_{c^0}$
such that $u y$ is of modified type
$(s_1)_c \cup (s_2)_{\bar{c}}$, we first
choose a number from each of the $r$ cycles of $u$ and
$r$ elements of $\G$, and then choose one additional new
number, say $i$, from one of the
$r$ cycles. (Keep in mind that we have counted every element twice
this way because we have made an
{\em ordered} choice of two $i$'s in
the same cycle of $x$.)
Then we arrange these $r+1$ numbers and
$r$ elements of $\G$ into a cycle
$y=( (h_{i_1}, \ldots, h_{i_r}, h_{i}), (i_1 \ldots i_r i))$
of modified type $(r)_{c^0}$
(we need $r+1$ elements of $\G$, however
the $i$-th element is uniquely
determined to be
$h_{i} =h_{i_1}^{-1} \ldots h_{i_r}^{-1}$
since the cycle-product is 1). In this way,
$uy$ will be of the modified type
of two cycles whose cycle-products
belong to opposite conjugacy classes,
by  Proposition~\ref{prop:constraint}~(2).

To ensure that $uy$ is of modified type
$(s_1)_c \cup (s_2)_{\bar{c}}$, we need to
choose the number $i$ such that
the two cycles on the right-hand side of
(\ref{eq:mult}) consist of an $s_1$-cycle and
an $s_2$-cycle. The number of choices of such
$i$ is $q(\mu,r, s)$ if $s_1 =s_2 =s$ and
 $2 q(\mu,r, s_1)$ if $s_1  \neq s_2. $
We can further write $ y =y_1 y_2$ where
$y_1 =( (h_{i_1}h_{i}, h_{i_2}, \ldots, h_{i_r}),
 (i_1 i_2 \ldots i_r))$ and
$ y_2 = ((h_{i}^{-1}, h_{i}), (i_r i)).$
According to Proposition~\ref{prop:constraint},
$u y_1$ is of the modified type $(n-1)_{c^0}$, and
$uy =(u y_1) y_2$ has cycle products lying in
opposite conjugacy classes $\{c, \bar{c}\}$ if and only
if $h_{i}=h_{i_1}^{-1} \ldots h_{i_r}^{-1}$ lies in $c
\text{ or } \bar{c}$.
Note that $\zeta_c = \zeta_{\bar{c}}.$

Taking all these considerations into account,
we obtain the number $A_3$ of choices for $y$
so that $uy$ belongs to $K_{(s_1)_c \cup (s_2)_{\bar{c}}}$:

 \begin{eqnarray*}
 A_3
  &=&   \left\{
      \everymath{\displaystyle}
      \begin{array}{ll}
 2^{m_1(\mu)} 3^{m_2(\mu)}
 \cdots  q(\mu,r;s)  |\G|^r/ (2 \zeta_c )
  , &\text{if }  c=\bar{c} \text{ and } s_1=s_2=s,  \\
  2^{m_1(\mu)} 3^{m_2(\mu)}
 \cdots  q(\mu,r;s_1)  |\G|^r/\zeta_c ,
   &  \text{otherwise}.
      \end{array}
    \right.
\end{eqnarray*}

The number $D$ of elements in $K_{(s_1)_c \cup (s_2)_{\bar{c}} }(n)$
is, by (\ref{eq:centord}), given by

 \begin{eqnarray*}
 D
  &=&   \left\{
      \everymath{\displaystyle}
      \begin{array}{ll}
{ |\G|^n n!}/(2 s^2 \zeta_c^2)
  , &\text{if }  c=\bar{c} \text{ and } s_1=s_2=s,  \\
 { |\G|^n n!} / (s_1s_2 \zeta_c^2),
  &  \text{otherwise}.
      \end{array}
    \right.
\end{eqnarray*}

Thus, $B =A_2 A_3/D$,
and the coefficient $N$ is given by
$N =A_1 B =A_1 A_2 A_3 /D.$
Remembering $n =\| \mu\| +r$ and
$q (\mu,r,s) =
{(r-\ell (\mu))! \prod_{i\ge 1}m_i(\mu)!} \cdot p(\mu,r,s)$,
we have established the theorem.
\end{proof}

\begin{example}  \label{ex:coeff}
 \begin{enumerate}
 \item
A specialization of Theorem~\ref{th:singlecyc2} gives us
$$ K_{(1)_{c^0}} K_{(1)_{c^0}}
 =2 K_{(1^2)_{c^0}} +3 K_{(2)_{c^0}}
 + \sum_{\{c,\bar{c}\},c \neq c^0} \zeta_c K_{(1)_{c} \cup (1)_{\bar{c}}}.
$$
 \item
 If we specialize Theorem~\ref{th:singlecyc1}
to $\la = (s)_{c^0}$, then we have two nontrivial
terms (i.e. $\mu =\emptyset$ or $\mu =(s)_{c^0}$):
 \begin{eqnarray*}
 K_{(s)_{c^0}} K_{(r)_{c}} =
 K_{(s)_{c^0}\cup (r)_{c}}
  + (s+r) K_{ (s +r)_{c}}, \qquad c \neq c^0.
 \end{eqnarray*}
By changing $r$ to $s$ in Theorem~\ref{th:singlecyc2}
and then specializing to $\la = (r)_c$, we obtain
the same answer (as expected).
\end{enumerate}
\end{example}

We now introduce a partial ordering on $\mathcal P_n (\G_*)$.
Given $\la, \mu \in \mathcal P_n (\G_*)$, $\mu \geq \la$ if
$\mu (c^0)$ is strictly contained in $\la (c^0)$,
or if $\mu (c) \geq \la (c)$ for each $c \in \G_*$
in the case when $\| \la (c) \| = \| \mu (c) \|$ for all $c$.
We can reformulate Theorem~\ref{th:singlecyc1} as follows.

\begin{theorem}  \label{th:reformulate1}
 Let $c \neq c^0$, $\la \in \mathcal P(\G_*)$, and
 $\widetilde{\la} =\la -\la (c^0)$. Then
 \begin{enumerate}
  \item
 if $\| \la (c^0)\| +r =k$, then
 \begin{eqnarray*}
  a_{\la\, (r)_c}^{\widetilde{\la} \cup (k)_c}
  &=&
  \left\{
      \everymath{\displaystyle}
      \begin{array}{ll}
   \frac{(m_k(\la(c))+1)\cdot k \cdot (r-1)!}{ (r-\ell(\la (c^0)))!
  \prod_{i \geq 1} m_i(\la (c^0))!}
  , &\text{if }  \ell(\la (c^0)) \leq r,  \\
   0, &  \text{otherwise}.
      \end{array}
    \right.
\end{eqnarray*}
  \item
 if $\| \la \| +r = \| \nu\|$, and if we write
$\nu = (\nu (c))_{c \in \G_*}$ and
the partition $\nu(c) =(\nu_1, \nu_2, \ldots)$, then
$$a_{\la (r)_c}^{\nu}
 = \sum a_{\mu\; (r)_c}^{(\nu_i)_c}
$$
summed over pairs $(i,\mu)$, where $\mu =\mu (c^0)$
and $\mu \cup \nu = \la \cup (\nu_i)_c.$
  \item
the coefficient
 $a_{\la (r)_c}^{\nu} = 0$ unless $\nu \geq \la \cup (r)_c$,
 and $a_{\la (r)_c}^{\la \cup (r)_c} >0$.
 \end{enumerate}
\end{theorem}

\begin{proof}
 The coefficient in part (1) is the one on the
right-hand side of the general formula in Theorem~\ref{th:singlecyc1}
corresponding to $\mu =\la(c^0)$. Part (2) is a simple
rewriting of this general formula. The first half of part (3)
follows from the definition of the partial order $\geq$.
Note that $a_{\la (r)_c}^{\la \cup (r)_c}$ is the
coefficient on the right-hand side of the general formula
corresponding to $\mu =\emptyset$.
\end{proof}

\begin{remark} \label{rem:order}
 It is possible also to reformulate Theorem~\ref{th:singlecyc2}
in a form like Theorem~\ref{th:reformulate1}. Let us simply
mention that
$a_{\la (r)_{c^0}}^{\nu} = 0$ unless $\nu \geq \la \cup (r)_{c^0}$,
 and $a_{\la (r)_{c^0}}^{\la \cup (r)_{c^0}} >0$.
\end{remark}

\begin{remark}
 In the symmetric group case, the product in the FH-ring
with a conjugacy class of single cycles was
computed in \cite{FH}. The wreath product
case here is considerably more difficult.
The reformulation in the form of Theorem~\ref{th:reformulate1}
is similar to the presentation of Macdonald
(cf. Ex.~24, \cite{Mac}, pp. 132).
Some other structure constants have also been computed
in the symmetric group case by Goulden-Jackson
and others (cf. \cite{GJ} and the references therein).
\end{remark}

The structure constants $a_{\la\mu}^\nu$ for $\nu = \la \cup \mu$
are easily computed just as in the symmetric group case
(cf. Lemma~3.10, \cite{FH}).

\begin{proposition}
 Let $\la =(r^{m_r(\la(c))})_{r \ge 1, c\in\G_*}$ and
$\mu =(r^{m_r(\mu(c))})_{r \ge 1, c\in\G_*}$. We have
$$ a_{\la\mu}^{\la \cup \mu}
 = \prod_{r \ge 1, c\in\G_*}
 {{m_r(\la(c)) + m_r(\mu(c))} \choose{m_r(\la(c))}}.
$$
\end{proposition}

\begin{proof}
Fix an element $z$ in the conjugacy class $\mathcal K_{\la \cup \mu}$.
The coefficient $ a_{\la\mu}^{\la \cup \mu}$ is equal to
the number of pairs $(x,y)$ such that $xy =z$, where $x \in \mathcal K_\la$
and $y \in \mathcal K_\mu$,
which is given by the formula in the proposition.
\end{proof}

%

%
%
%
\subsection{The structures of the Farahat-Higman ring}

\begin{theorem}  \label{th:free}
  Let $\la = (\la(c))_{c\in \G_*}$
 be in $\mathcal P (\G_*)$, where
 $\la(c) =(\la_1(c), \la_2(c), \ldots).$ Then the product
 $\prod_{c\in \G_*} K_{\la_1}(c) K_{\la_2}(c) \ldots$
 in the FH-ring $\mathcal G_\G$ is of the form
 $$\prod_{c\in \G_*} K_{\la_1}(c) K_{\la_2}(c) \ldots
 =\sum_{\mu \geq \la} d_{\la \mu} K_\mu$$
with $d_{\la \mu} \in \Z_+$
and $d_{\la\la} >0$. In particular, $K_r(c)$,
where $r \ge 1, c \in \G_*$, are algebraically
independent elements of the FH-ring $\mathcal G_\G$,
and generate $\mathbb Q \otimes_\Z \mathcal G_\G$.
\end{theorem}

\begin{proof}
 The first statement follows from Theorems~\ref{th:reformulate1}~(3)
 and Remark~\ref{rem:order}. The second statement
holds since the matrix $(d_{\la\mu})$ is triangular
and nonsingular.
\end{proof}

\begin{remark}
Theorems~\ref{th:singlecyc1}
and \ref{th:singlecyc2} also provide a recursive way
of computing the general structure constants $a_{\la\mu}^\nu$.
\end{remark}

Denote by $\mathcal I_\G(n)$ the subspace of
$\mathcal G_\G(n)$ which is spanned by those
$K_\la(n)$ where $\la \neq \la(c^0)$. Similarly,
we introduce the subspace $\mathcal I_\G$
of $\mathcal G_\G$. Recall that $\mathcal G(n)$ and $\mathcal G$
denote the rings $\mathcal G_\G(n)$
and $\mathcal G_\G$ when $\G$ is trivial.

\begin{corollary} \label{cor:quotient}
The space $\mathcal I_\G(n)$ is a graded ideal of the ring
$\mathcal G_\G(n)$, and the quotient ring
$\mathcal G_\G(n) / \mathcal I_\G(n)$
is canonically isomorphic to the ring $\mathcal G(n)$.
Similarly,  $\mathcal I_\G$ is a graded ideal of the ring
$\mathcal G_\G$, and the quotient ring
$\mathcal G_\G / \mathcal I_\G$
is canonically isomorphic to the ring $\mathcal G$.
\end{corollary}

\begin{proof}
It suffices to prove the statements for $\mathcal G_\G$
as the same proof applies to $\mathcal G_\G(n)$.

We proceed by induction on the length of elements
of the set $\mathcal P(\G_*)$. By observation from the formulas in
Theorems~\ref{th:singlecyc1} and \ref{th:singlecyc2},
we see that $\mathcal I_\G$ is closed
under the multiplication by any conjugacy class $K_\rho$
of single-cycle type (i.e. $\ell (\rho) =1$).

For a given $\rho \in \mathcal P(\G_*)$ with $\ell (\rho) >1$,
let us assume that $\mathcal I_\G$ is closed
under the multiplication by $K_\la$ for all $\la$ with
$\ell(\la) < \ell (\rho)$.

First assume that some
cycle-product of $\rho$ is not $1$, i.e., $\rho \neq \rho (c^0)$.
Then we can write $\rho = \la \cup (r)_c$
for some $c \neq c^0$ and  some $\la$ with $\ell(\la) <\ell(\rho)$.
By applying Theorem~\ref{th:singlecyc1},
we observe that all the conjugacy classes
appearing in the product $K_\la K_{(r)_{c}}$
have length less than $\ell (\rho) =\ell (\la) +1$ except
for the term $K_{\la \cup (r)_c} =K_\rho.$
By applying the induction assumption to
$K_\la, K_{(r)_c}$, and all those conjugacy classes in
the product $K_\la K_{(r)_{c}}$ except for $K_\rho$, we see that
$\mathcal I_\G$ is closed under the
multiplication by $K_\rho$.

Now assume that $\rho =\rho (c^0)$. We can write
$\rho = \la \cup (r)_{c^0}$
for some $\la =\la(c^0)$ with $\ell(\la) <\ell(\rho)$.
By applying Theorem~\ref{th:singlecyc2},
we observe that all the conjugacy classes
appearing in the first term on the right-hand
side of the product $K_\la K_{(r)_{c^0}}$
have length less than $\ell (\rho) =\ell (\la) +1$ except
for the term $K_{\la \cup (r)_{c^0}} =K_\rho;$
the second term is $0$ since $\la =\la(c^0)$; we
further observe that  all the conjugacy classes
appearing in the remaining three terms
on the right-hand side of the product $K_\la K_{(r)_{c^0}}$
have length less than or equal to $\ell (\rho)$,
and each of the corresponding modified types contains some
cycle-product which is not $1$. Thus we can
apply the induction
assumption together with the case established
in the previous paragraph. This proves that
$\mathcal I_\G$ is closed under the
multiplication by $K_\rho$.

This finishes the proof that
$\mathcal I_\G$ is a graded ideal of the ring $\mathcal G_\G$.
It now follows from Theorem~\ref{th:free}
and the explicit formula in Theorem~\ref{th:singlecyc2}
that the obvious map between the quotient ring
$\mathcal G_\G / \mathcal I_\G$ and the ring $\mathcal G$
is a ring isomorphism.
\end{proof}

\begin{remark}
 Another proof of Corollary~\ref{cor:quotient} goes as follows.
Regarding $S_n$ as a subgroup of $\Gn$, we note that
the intersection $\mathcal K_\la(n) \cap S_n$ of $S_n$ with
a conjugacy class $\mathcal K_\la(n)$ of $\Gn$ is a conjugacy class of $S_n$
when $\la =\la (c^0)$, and is $\emptyset$ otherwise.
Thus we have a surjective ring homomorphism
$\varphi_n: R_\Z(\Gn) \rightarrow R_\Z(S_n)$, by
$K_\la(n) \mapsto K_\la(n) \cap S_n$ for $\la =\la (c^0)$
and $K_\la(n) \mapsto 0$ otherwise.
The ring homomorphism $\varphi_n$ is compatible with
the filtrations and thus induces a surjective ring homomorphism
$\phi_n :\mathcal G_\G(n) \rightarrow \mathcal G (n).$
It follows that the kernel of $\phi_n$ coincides with $\mathcal I_\G(n)$
and thus we have the ring isomorphism
$\mathcal G_\G(n) /\mathcal I_\G(n) \cong\mathcal G (n)$.
\end{remark}

One may look at the Example~\ref{ex:coeff}, where
the type of phenomena in Corollary~\ref{cor:quotient} is manifest.

The next corollary can be established by the same type of induction
argument as in the first proof of Corollary~\ref{cor:quotient},
which we omit here.

\begin{corollary}  \label{cor:generator}
 The $n|\G_*|$ elements $K_{(r)_c}(n)$, where $0 \leq r <n$
for $c =c^0$ and $1 \leq r \leq n$
for $c \neq c^0$, form  a set of ring generators of
the ring $\mathbb Q \otimes_\Z \mathcal G_\G(n)$.
\end{corollary}

\begin{remark}
Corollary~\ref{cor:generator} implies that
the same $n|\G_*|$ elements also form
a set of ring generators of the ring $R(\Gn)$. This latter
fact was established in \cite{Wa4}, Theorem~5.9~(ii)
in a totally different approach.
\end{remark}

\begin{definition}
A conjugacy class $c \in \G_*$ is called {\em real} (resp. {\em complex})
if $c = \bar{c}$ (resp. $c \neq \bar{c}$).
Define $\G^{cx}_* =\{c \in \G_* \mid c \neq \bar{c} \}$
and  $\G^{re}_* =\{c \in \G_* \mid c \neq c^0 \text{ and } c = \bar{c} \}.$
\end{definition}

Of course, we have $|\G_*| = |\G^{cx}_*| + |\G^{re}_*| +1$ and
$|\G^{cx}_*|$ is always an even integer.

\begin{theorem}  \label{th:numinv}
The ring $\C \otimes_\Z \mathcal G_\G(n)$ for every $n$ does not depend
on the finite group $\G$, but depends only on
the two numbers $|\G_*|$ and $|\G^{re}_*|$.
In particular, if $G$ is another finite group
such that $|\G_*| = |G_*|$ and $|\G^{re}_*| = |G^{re}_*|$,
then the rings $\C \otimes_\Z \mathcal G_\G(n)$ and
$\C \otimes_\Z \mathcal G_G (n)$ are isomorphic for every $n$.
\end{theorem}

\begin{proof}
We observe that the factor $\zeta_c$ in the formula in
Theorem~\ref{th:singlecyc2} can be annihilated by a rescaling
$\widetilde{K}_\la = \prod_{c \neq c^0}
 \zeta_c^{\ell(\la(c))/2} K_\la$ over $\C$.
It follows from Theorems~\ref{th:singlecyc1} and \ref{th:singlecyc2} that
the ring $\C \otimes_\Z \mathcal G_\G(n)$ with respect
to the new basis $\widetilde{K}_\la$
no longer depends on the group $\G$.

Let $G$ be another finite group
such that $|\G_*| = |G_*|$ and $|\G^{re}_*| = |G^{re}_*|$.
Fix a bijection $\varphi$ from $\G_*$ to $G_*$ such that
its restriction to $\G^{re}_*$ is one-to-one and onto
$G^{re}_*$. Then the obvious map by matching
the rescaled conjugacy classes of $\Gn$ and
those of $G_n$ in terms of  $\varphi$
is a ring isomorphism by Theorems~\ref{th:singlecyc1},
\ref{th:singlecyc2} and Corollary~\ref{cor:generator}.
\end{proof}

The following conjecture concerns about
the converse of Theorem~\ref{th:numinv}.

\begin{conjecture}   \label{conj:twonumber}
Let $\G$ and $G$ be two finite groups.
If the rings $\C \otimes_\Z \mathcal G_\G(n)$ and
$\C \otimes_\Z \mathcal G_G (n)$ are isomorphic for every $n$,
then $|\G^{re}_*| = |G^{re}_*|$, (in addition to the obvious identity
$|\G_*| = |G_*|$).
\end{conjecture}
\subsection{Real conjugacy classes and the dual McKay correspondence}

The classification of finite subgroups of $SL_2(\C)$
is well known. The following is a complete list:
cyclic groups $\Z_n$ of order $n$,
the binary dihedral groups $BT_{4n}$ of order $4n$ $(n \geq 3)$, the binary
tetrahedral group $BT$, the binary octahedral group $BO$, and
the binary icosahedral group $BI$. It is well known
(cf. e.g. \cite{Cox, Slo}) that
there is a bijection between finite subgroups $\G$ of $SL_2(\C)$
(modulo conjugation) and Dynkin diagrams $\Delta$ of ADE types.

A nice presentation of these finite groups in
terms of generators and relations can be found
in Coxeter \cite{Cox}. In the table I below,
we have set $Z=-I$, the negative of the identity matrix in $SL_2(\C)$.
So, $Z^2=1$. The cyclic group has one generator $A$,
the binary dihedral group has two generators $A,B$, and
each of the binary polyhedral groups has three generators $A,B,C$.

\vskip 1.5pc
\centerline{{\bf Table I}}
\vskip 1pc

\vbox{\tabskip=0pt \offinterlineskip
\def\tablerule{\noalign{\hrule}}
\halign to385pt{\strut#& \vrule#
\tabskip=1em plus2em&
\hfil#& \vrule#&\hfil#\hfil& \vrule#& #\hfil& \vrule#&
#\hfil& \vrule#\tabskip=0pt \cr\tablerule
&&\omit\hidewidth $\Delta$\hidewidth&&
\omit\hidewidth $|\G|$\hidewidth&&
\omit\hidewidth $\G$\hidewidth&&
\omit\hidewidth Relations\hidewidth&\cr\tablerule

&&$A_{n-1}$&& $n$&&cyclic $\Z_n$ &&$A^n =1$&\cr\tablerule
&&$D_{n+2}$&&$4n$&&binary dihedral $BD_{4n}$
&&$A^n =B^2=(AB)^2=Z$&\cr\tablerule
&&$E_6$&&$24$&&binary tetrahedral $BT$ &&$A^3 =B^3=C^2=Z$&\cr\tablerule
&&$E_7$&&$48$&&binary octahedral $BO$ &&$A^4 =B^3=C^2=Z$&\cr\tablerule
&&$E_8$&&$120$&&binary icosahedral $BI$ &&$A^5 =B^3=C^2=Z$&\cr\tablerule}}
\vskip 2pc

There is a dual McKay correspondence
introduced by Ito-Reid \cite{IR} which identify
bijectively the nontrivial conjugacy classes of $\G$ and
the vertices of the corresponding Dynkin diagram.
Brylinski \cite{Bry} later gave a more transparent and canonical
construction, by using Mumford's description \cite{Mum} of $\G$
as the fundamental group of the complement of the exceptional divisor
of the minimal resolution of $\C^2/\G$. Following \cite{Bry}, we have
attached the specific representatives of
conjugacy classes to the vertices as follows.

 \begin{equation*}
 \begin{picture}(150,75) 
 \put(-90,18){$A_{n-1}$}
 \put(-30,20){$\bigcirc$}
 \put(-16,23){\line(1,0){23}}
 \put(10,20){$\bigcirc$}
 \put(24,23){\line(1,0){15}}
 \put(41,22){ \dots }
 \put(64,23){\line(1,0){15}}
 \put(82,20){$\bigcirc$}
 \put(-30,5){$A$}
 \put(12,5){$A^2$}
 \put(80,5){$A^{n-1}$}
 \end{picture}
 \end{equation*}
%
%
 \begin{equation*}
 \begin{picture}(150,75) 
 \put(-90,28){$D_{n+2}$}
 \put(-40,30){$\bigcirc$}
 \put(-26,33){\line(1,0){23}}
 \put(0,30){$\bigcirc$}
 \put(14,33){\line(1,0){15}}
 \put(35,30){$\cdots$}
 \put(52,33){\line(1,0){15}}
 \put(68,32){$\bigcirc$ }
 \put(81,33){\line(1,0){23}}
 \put(105,30){$\bigcirc$}
 \put(116,38){\line(1,1){18}}
 \put(135,60){$\bigcirc$}
 \put(116,27){\line(1,-1){18}}
 \put(135,0){$\bigcirc$}
 \put(-40,15){$A$}
 \put(0,15){$A^2$}
 \put(72,15){$A^{n-1}$}
 \put(106,15){$Z$}
 \put(148,0){$B$}
 \put(148,60){$BA$}
 \end{picture}
 \end{equation*}
%
%
 \begin{equation*}
 \begin{picture}(150,75) 
 \put(-90,28){$E_6$}
 \put(-40,30){$\bigcirc$}
 \put(-26,33){\line(1,0){23}}
 \put(0,30){$\bigcirc$}
 \put(14,33){\line(1,0){23}}
 \put(40,68){$\bigcirc$}
 \put(46,41){\line(0,1){23}}
 \put(40,30){$\bigcirc$}
 \put(54,33){\line(1,0){23}}
 \put(80,30){$\bigcirc$}
 \put(94,33){\line(1,0){23}}
 \put(120,30){$\bigcirc$}
 \put(-40,15){$B$}
 \put(2,15){$B^2$}
 \put(43,15){$Z$}
 \put(82,15){$A^2$}
 \put(123,15){$A$}
 \put(54,67){$C$}
 \end{picture}
 \end{equation*}
%
%
 \begin{equation*}
 \begin{picture}(150,75) 
 \put(-90,28){$E_7$}
 \put(-40,30){$\bigcirc$}
 \put(-26,33){\line(1,0){23}}
 \put(0,30){$\bigcirc$}
 \put(14,33){\line(1,0){23}}
 \put(40,68){$\bigcirc$}
 \put(46,41){\line(0,1){23}}
 \put(40,30){$\bigcirc$}
 \put(54,33){\line(1,0){23}}
 \put(80,30){$\bigcirc$}
 \put(94,33){\line(1,0){23}}
 \put(120,30){$\bigcirc$}
 \put(134,33){\line(1,0){23}}
 \put(160,30){$\bigcirc$}
 \put(-40,15){$B$}
 \put(2,15){$B^2$}
 \put(43,15){$Z$}
 \put(82,15){$A^3$}
 \put(123,15){$A^2$}
 \put(164,15){$A$}
 \put(54,67){$C$}
 \end{picture}
 \end{equation*}
%
%
 \begin{equation*}
 \begin{picture}(150,75) 
 \put(-90,28){$E_8$}
 \put(-40,30){$\bigcirc$}
 \put(-26,33){\line(1,0){23}}
 \put(0,30){$\bigcirc$}
 \put(14,33){\line(1,0){23}}
 \put(40,68){$\bigcirc$}
 \put(46,41){\line(0,1){23}}
 \put(40,30){$\bigcirc$}
 \put(54,33){\line(1,0){23}}
 \put(80,30){$\bigcirc$}
 \put(94,33){\line(1,0){23}}
 \put(120,30){$\bigcirc$}
 \put(134,33){\line(1,0){23}}
 \put(160,30){$\bigcirc$}
 \put(174,33){\line(1,0){23}}
 \put(200,30){$\bigcirc$}
 \put(-40,15){$B$}
 \put(2,15){$B^2$}
 \put(43,15){$Z$}
 \put(82,15){$A^4$}
 \put(123,15){$A^3$}
 \put(164,15){$A^2$}
 \put(205,15){$A$}
 \put(54,67){$C$}
 \end{picture}
 \end{equation*}

\begin{theorem}  \label{th:selfdual}
\begin{enumerate}
\item
 Via the dual McKay correspondence,
the subset $\G^{re}_*$ (resp. $\G^{cx}_*$) of conjugacy classes
for $\G \leq SL_2(\C)$ corresponds to
the set of vertices in the Dynkin diagrams
which are fixed (resp. not fixed) by the diagram automorphism $\tau^\G$.
Here $\tau^\G$ is the automorphism which exchanges the two end-points of
Dynkin diagrams for $\G =\Z_n$ $(n \geq 2)$, $BD_{8k+4}$, and $BT$,
and is the trivial automorphism otherwise.

\item
The set $\G^{cx}_*$ is given as follows:
$\{\langle A^{i}\rangle, 1 \leq i <2k+1 \}$ for $\Z_{2k+1}$ $(k>0)$;
$\{\langle A^{i}\rangle, 1 \leq i <2k, i \neq k \}$ for $\Z_{2k}$;
$\{\langle B\rangle, \langle BA \rangle\}$ for $BD_{8k+4}$ $(k>0)$;
$\{\langle A\rangle,\langle B\rangle,
\langle A^2\rangle, \langle B^2\rangle\}$ for $BT$;
the empty set $\emptyset$ for $BD_{8k}$, $BO$ and $BI$.

\item
The cardinalities of
$\G_*$, $\G^{re}_*$, $\G^{cx}_*$, together with the order
of $\tau^\G$, are listed in the table II below.

\vskip 1.5pc
\centerline{{\bf Table II}}
\vskip 1pc

\vbox{
\vspace{.2cm}
\tabskip=0pt \offinterlineskip
\def\tablerule{\noalign{\hrule}}
\halign to395pt{\strut#& \vrule#
\tabskip=1em plus2em&
\hfil#& \vrule#&\hfil#\hfil& \vrule#& #\hfil& \vrule#
&#\hfil& \vrule#
&#\hfil& \vrule#
&#\hfil& \vrule#
&#\hfil& \vrule#
&#\hfil& \vrule#\tabskip=0pt \cr\tablerule
&&\omit\hidewidth $\G$\hidewidth&&
\omit\hidewidth $\Z_{2k+1}$ \hidewidth&&
\omit\hidewidth $\Z_{2k}$ \hidewidth&&
\omit\hidewidth $BD_{8k}$\hidewidth&&
\omit\hidewidth $BD_{8k+4}$ \hidewidth&&
\omit\hidewidth $BT$  \hidewidth&&
\omit\hidewidth $BO$ \hidewidth&&
\omit\hidewidth $BI$ \hidewidth&\cr\tablerule

&&$|\G_*|$  &&$2k+1$&&$2k$  &&$2k+3$&&$2k+4$&& $7$&& $8$ && $9$ &\cr\tablerule
&&$|\G^{cx}_*|$&&$2k$&&$2k-2$&& $0$ && $2$&& $4$&& $0$ && $0$ &\cr\tablerule
&&$|\G^{re}_*|$&& $0$ &&$ 1 $ &&$2k+2$&&$2k+1$&& $2$&& $7$&& $8$&\cr\tablerule
&&$o(\tau^\G)$&& $2$ && $2$ &&$1$ &&$2$ && $2$ && $1$ && $1$ &\cr\tablerule
}}
\end{enumerate}
\end{theorem}

\begin{proof}
Parts (1) and (3) easily follow from (2) by
comparing with with dual McKay correspondence described
above. So let us prove (2) below.
The cyclic group case is evident.

For the groups
$BD_{4n}$, one can easily write down all $n+3$ conjugacy
classes: $\{ 1\}$, $\{ Z\}$,
$\{A^{\pm i} \}$ $(1\leq i \leq n-1)$,
$\{A^{2j}B, 0\leq j \leq n-1 \}$, and
$\{A^{2j-1}B, 0\leq j \leq n-1 \}$.
The results for $BD_{4n}$ (depending on whether
$n$ is even or odd) follow now
from this description together with the fact
that $B^{-1} =A^nB$.

For the binary tetrahedral group $BT$, we use
the explicit elements given in terms of quaternions in
Coxeter \cite{Cox}, (7.26), pp. 77.
They are $\pm 1, \pm i, \pm j, \pm k,$ and
$( \pm 1\pm i\pm j \pm k)/2$, where $i,j,k$ stands for
the quaternion generators. In terms of the notations
in Table I, we have $A =(1+i+j+k)/2$, $C=i$ and $B =(1+i-j+k)/2$
(cf. (7.22) in \cite{Cox}).
Then a direct calculation
leads to the following list of 7 conjugacy classes of $BT$:
$\{1 \}, \{ -1 \}, K_0, K_1, K_2,- K_1,$ and $ -K_2$,
where
\begin{eqnarray*}
K_0 &=& \{ \pm i, \pm j, \pm k \} \\
K_1 &=& \{ (1+i+j+k)/2, (1+i-j-k)/2,   (1-i+j-k)/2,
 (1-i-j+k)/2 \} \\
K_2 &=& \{ (1-i-j-k)/2, (1-i+j+k)/2,   (1+i-j+k)/2,
 (1+i+j-k)/2 \}
\end{eqnarray*}
and $- K_i$
consists of the negatives of all elements in $K_i$, $i=1,2$.
One verifies that
$\langle A\rangle =K_1, \langle A^{-1}\rangle =\langle B\rangle =K_2$,
and similarly, $\langle A^2 \rangle = -K_2$ and
$\langle A^{-2}\rangle =\langle B^2\rangle = -K_1$.
Thus, we have $(BT)_*^{cx} =\{\langle A\rangle,\langle B\rangle,
\langle A^2\rangle, \langle B^2\rangle\}$.

For the binary icosahedral group $BI$, we note that the orders of
the representatives $A, A^2, A^3, A^4, B, B^2$ and $C$
of the conjugacy classes
as marked in the $E_8$ Dynkin diagram
are respectively $10, 5, 10, 5, 6, 3$ and $4$.
To show that $(BI)^{cx}_* =\emptyset$,
it remains to check that the elements of the same order
do not lie in opposite conjugacy classes.
Let us check that $A$ is not conjugate to
the inverse of $A^3$. If $A$ were conjugate to $A^{-3}$, then
$A^3$  is conjugate to $(A^{-3})^3 =A$ (since $A^{10}=1$)
which contradicts
with the fact that $A$ and $A^3$ are representatives
of different conjugacy classes. Similarly,
we see that $A^2 $ is not conjugate to the inverse of $ A^4$.

This argument also applies to the binary
octahedral group $BO$. Alternatively,
one can also verify the theorem for $BO$ directly using
its explicit realization in quaternions
given by (7.26) and (7.35) in Coxeter \cite{Cox}.
\end{proof}

\begin{corollary}  \label{cor:amusing}
 The rings $\C \otimes_\Z \mathcal G_{BD_{24}}(n)$
and $\C \otimes_\Z \mathcal G_{BI}(n)$
are isomorphic for each $n$.
\end{corollary}

\begin{proof}
Follows from Theorem~\ref{th:numinv} and Theorem~\ref{th:selfdual}.
\end{proof}

The group $\G$ acts on $\C^2 -\{0\}$ freely and $\C^2/\G$
has a simple singularity at the origin. Let us denote
by $\ale$ the minimal resolution of $\C^2/\G$.
The bijection between the finite subgroups of $SL_2(\C)$
and the Dynkin diagrams of ADE types can also be
seen through the geometry of minimal resolutions.
A classical result of Du Val (cf. e.g. \cite{Slo})
says that the set of
exceptional curves in the minimal resolution
$\ale$ corresponds naturally to the set
of vertices of the corresponding Dynkin diagram.

\begin{remark}  \label{rem:autom}
Combining with the dual McKay correspondence \cite{IR, Bry}, we
may ask what is the geometric significance of the set of
exceptional curves corresponding to the set $\G^{re}_*$, or
alternatively what is the geometric significance of the Dynkin
diagram automorphisms $\tau^\G$ in Theorem~\ref{th:selfdual}. It
turns out that this question has the following beautiful answer
which was suggested by Dolgachev to the author. Consider the
automorphism $\tau_\G$ of $\ale$ which comes from the one which
sends $z$ to $ -z$ if we realize the simple singularities as the
hypersurface singularities in $\C^3$ in the standard way (cf. e.g.
Slodowy \cite{Slo}): $x^n +y^2+z^2 =0$ for $\G =\Z_n$;
$x(y^2-x^n)+z^2 =0$ for $\G =BD_{4n}$; $x^4 +y^3 +z^2 =0$ for $\G
=BT$; $x^3 +xy^3 +z^2 =0$ for $\G =BO$; and $x^5 +y^3 +z^2 =0$ for
$\G=BI$. Then $\tau_\G$ induces an automorphism on the degree-$2$
homology group of $\ale$ which fixes exactly the homology classes
associated to the exceptional curves corresponding to $\G_*^{re}$,
i.e. this induced automorphism is exactly the $\tau^\G$ in
Theorem~\ref{th:selfdual}. In addition, $\tau_\G$ is an involution
of $\C^2/\G$ such that $\tau_\G$ and $\G$ generate a complex
reflection subgroup $\widehat{\G}$ of $GL_2(\C)$. The quotient
$\C^2/\widehat{\G}$ is isomorphic to the affine plane.

On the other hand, it is brought to my attention by Brian Parshall
that the diagram automorphism $\tau^\G$ can be interpreted as
minus the longest reflection in the corresponding Weyl groups
acting on the set of simple roots, if we attach a simple root to
each of the vertices of the Dynkin diagrams associates to $\G$.
\end{remark}

\section{Connections with the cohomology rings of
Hilbert schemes}
\label{sec:scheme}
\subsection{Minimal resolutions and Hilbert schemes}

Given a quasi-projective surface $X$, we denote by
$\Xn$ the Hilbert scheme of $n$ points on $X$.
An element in $\Xn$ is represented by a length $n$
zero-dimensional closed subscheme. A
well-known theorem of Fogarty states
that $\Xn$ is nonsingular and the Hilbert-Chow
morphism $\pi_n: \Xn \longrightarrow X^n/S_n$, which
sends an element to its support, is a resolution
of singularities (cf. \cite{Na3}).

Connections between Hilbert
schemes and wreath products were first pointed out
in \cite{Wa1}. In the special case when $X$
is the minimal resolution $\ale$ of a simple
singularity $\C^2/\G$ (for $\G \leq SL_2(\C)$), we have
the following commutative diagram:

$$\CD (\ale)^{[n]} @>>> (\ale)^n  /S_n \\ @VVV @VVV \\ \C^{2n}  /\Gn @<{\cong}<<
(\C^2/\G)^n /S_n \endCD $$
which defines a resolution
of singularities $(\ale)^{[n]} \rightarrow \C^{2n}  /\Gn.$
It has thus been expected that all the geometric invariants
of $(\ale)^{[n]}$
can be described entirely by using the wreath product $\Gn$.

\subsection{Shift numbers for the wreath product orbifolds}

Given a complex manifold $M$ acted upon by
a finite group $G$, Zaslow \cite{Zas} introduced
a shift number associated to $g\in G$ (depending on
the connected components of the fixed-point set $M^g$).
In the case when $M$ is the affine space $\C^N$,
also see Ito-Reid \cite{IR}. For the orbifold $\C^{2n}/\Gn$,
the shift numbers have been computed
in \cite{WaZ}. Our observation is the following.

\begin{lemma}  \label{lem:samenumb}
Let $\G$ be a finite subgroup of $SL_2(\C)$.
Let $g \in \Gn$ be of modified type $\rho\in \mathcal P(\G_*)$.
The shift number for $g$ in the orbifold $\C^{2n}/\Gn$
coincides with $\Vert \rho\Vert$.
\end{lemma}

\begin{proof}
 The shift numbers for a general wreath product orbifold were
calculated in formula (9), pp. 162, \cite{WaZ}.
Let us pinpoint what the notations in {\em loc. cit.} mean, when
we specialize to the affine orbifold $\C^{2n}/\Gn$:
$d=2$ (which is the dimension of $\C^2$);
the fixed-point set of an element
$a \in \G$ acting on $\C^2$ is a point
(resp. the whole $\C^2$) and
thus the shift number in $\C^2/\G$ is $F^c=1$ (resp. $F^c=0$)
if $a \neq 1 $ (resp. $a = 1 $); in particular, the fixed-point set
is always connected which means $N_c=1$.
This shows that the formula (9) in {\em loc. cit.}
agrees exactly with our definition of $\Vert \rho\Vert$.
\end{proof}

\begin{remark}
A general wreath product $\Gn$, associated to
a finite subgroup of $\G \leq SL_k(\C)$, acts on $\C^{kn}$.
One can still define the shift numbers in the sense of
Zaslow. But they in general do not coincide with
(a constant multiple of) the values of $\| \cdot \|$
except some rare cases. Such cases
include $V^n/S_n$ or  $V^n/ \Gn$ (where
$\G =\Z/2\Z$ acts as $\pm \text{Id}_V$),
for a vector space $V =\C^{2m}$, $m \ge 2$.
It will be very interesting to analyze the
structure of the graded algebra of
$R_\Z(\Gn)$ associated to this orbifold filtration
for a general $\G$.
\end{remark}
\subsection{A ring isomorphism}

Lehn-Sorger \cite{LS1} and independently
Vasserot \cite{Vas} gave a description of
the cohomology ring $H^*(\Xn)$ for $X =\C^2$
in terms of the symmetric group $S_n$.
(However, the reference \cite{FH} has never featured
in the literature on Hilbert schemes.)
In light of Lemma~\ref{lem:samenumb},
a theorem of Etingof-Ginzburg \cite{EG}, when
combined with the results of Nakajima etc, leads to the following
description of
the cohomology ring $H^*((\ale)^{[n]})$ with $\C$-coefficient
(cf. the remarks in a footnote  in the Introduction).

\begin{theorem}  \label{th:coh}
Let $\G$ be a finite subgroup of $SL_2(\C)$. There exists
a ring isomorphism between
the cohomology ring $H^*((\ale)^{[n]})$ with $\C$-coefficient
and the ring $\C \otimes_\Z \mathcal G_\G (n)$.
\end{theorem}

Combined with Theorem~\ref{th:coh},
the results in the previous sections provide
finer structures on the cohomology ring of $(\ale)^{[n]}$.
Let us restrict ourselves to mention the following corollaries.

\begin{corollary}   \label{cor:inv}
For each $n$, the cohomology rings of the Hilbert schemes of
$n$ points of the minimal resolution $\ale$
are determined by the numbers $|\G_*|$ and $|\G^{re}_*|$.
\end{corollary}

\begin{proof}
Follows from Theorem~\ref{th:coh}
and Theorem~\ref{th:numinv}.
\end{proof}

\begin{remark}
A theorem of \cite{LQW2} says that the cohomology
ring of $\Xn$ associated to a {\em projective}
surface $X$ depends only on the cohomology ring
of the surface $X$ and the canonical class $K_X$.
A natural and important question is what invariants on a
quasi-projective surface $X$ determine
the cohomology ring of $\Xn$ for all $n$.
Corollary~\ref{cor:inv} can be regarded as a first step in this direction.
\end{remark}

\begin{corollary}
The cohomology rings of the Hilbert schemes of
$n$ points of the minimal resolutions
for the binary group $BD_{24}$ of order $24$
and the binary Icosahedral group $BI$
are isomorphic for all $n$.
\end{corollary}

\begin{proof}
Follows from Theorem~\ref{th:coh}
and Corollary~\ref{cor:amusing}.
\end{proof}

Conjecture~\ref{conj:twonumber}, if true, would imply that
$BD_{24}$ and $BI$ are the only pair of finite subgroups
of $SL_2(\C)$ whose corresponding Hilbert schemes of $n$ points
have isomorphic cohomology ring structure for all $n$.
\subsection{Implications on the cohomology rings of Hilbert schemes}
\label{sec:predict}

In recent years, there has been much progress on the understanding
of the cohomology rings of Hilbert schemes $\Xn$ of $n$ points on
a (quasi-)projective surface $X$. Many results have been obtained
for a general projective surface $X$ using a vertex operator
approach \cite{Lehn, LQW1, LS2, LQW2, LQW3} built on \cite{Na2,
Na3, Gro}. However, for $X$ quasi-projective, the vertex operator
techniques are not directly applicable, and little is known about
the ring $H^*(\Xn)$, with the notable exceptions when $X$ is the
affine plane $\C^2$ or $\ale$, cf. \cite{LS1, Vas} and
Theorem~\ref{th:coh}.

As suggested in \cite{Wa3}, the structures of the class algebras
of wreath products $\Gn$ associated to an arbitrary finite group
$\G$ (resp. the graded ring associated to some appropriate
filtration) very much reflect those of the cohomology rings of
Hilbert schemes $\Xn$ of points on an arbitrary surface $X$ which
is projective (resp. quasi-projective). According to this
philosophy, many of the results obtained on the rings $\mathcal
G_\G(n)$ (and $\mathcal G_\G$) are expected to find their
counterparts (some of which could be quite subtle) on the
cohomology rings $H^*(\Xn)$ for a simply-connected
quasi-projective surface $X$. As an example for illustration, we
will formulate a conjecture below.

Recall that a Heisenberg algebra was constructed \cite{Na2, Na3,
Gro} geometrically which acts on $\oplus_{n\geq 0} H^*(\Xn)$
irreducibly with the vacuum vector $|0\rangle =1 \in
H^*(X^{[0]})$. For the sake of notational simplicity, let us
assume that $X$ is simply-connected and so it has no odd degree
cohomology classes. For $X$ either projective or quasi-projective,
the half of the Heisenberg algebra, i.e. the creation operators,
which are the only part we need below can always be modelled on
the ordinary cohomology group of $X$ (compare \cite{Na3}; cf.
\cite{Lehn}, Sect.~4.4), and they will be denoted by $\mathfrak
a_{-n}(\alpha)$, where $n \in \mathbb N, \alpha \in H^*(X)$.
Let us fix a linear basis $B_2$ of $H^2(X)$, so a linear
basis $S$ of $H^*(X)$ is given by $S =\{1_X\} \cup B_2$
for $X$ quasi-projective (and with the extra class
$\{[pt]\}$ for $X$ projective), where $1_X$ denotes
the unit. This leads to a natural linear basis
for $H^*(\Xn)$ (cf. \cite{Na3}), which we will refer to
as the Heisenberg monomial basis.

Now let $X$ be in addition quasi-projective. Just as the conjugacy
classes of $\Gn$, the Heisenberg monomial basis for $H^*(\Xn)$ can
be also parametrized by modified types as follows. Given $\mu =
(\mu (c))_{c \in S} \in \mathcal P(S)$, we denote by $\tilde{\mu}
= (\tilde{\mu} (c))_{c \in S} \in \mathcal P(S)$ with $\tilde{\mu}
(c) = \mu (c)$ for $c \in B_2$, and $\tilde{\mu} (1_X)
=(r^{m_r(\tilde{\mu}(1_X))})_{r\geq 1}$, where
$m_r(\tilde{\mu}(1_X)) =m_{r-1} ({\mu}(1_X))$ for $r \geq 2$ and
$m_1(\tilde{\mu}(1_X)) =n -\|\mu\|- \ell(\mu)$. (Here we have
denoted $ {\mu}(1_X) = (r^{m_r({\mu}(1_X))})_{r\geq 1}$).

We denote by
\begin{eqnarray*}
 \mathfrak h_\mu(n) =
 \frac1{(n-\|\mu\|- \ell(\mu))!}
  \prod_{c \in S}\prod_{r \ge 1}
  {\mathfrak a}_{-r}(c)^{m_r(\tilde{\mu}(c) )} \,|0\rangle
\end{eqnarray*}
if $n \geq \|\mu\| +\ell(\mu)$, and $0$ otherwise.
Then, those nonzero $\mathfrak h_\mu(n)$'s, as $\mu$
runs over $\mathcal P_n(S)$, form a linear
basis for $H^*(\Xn)$. (The $\mathfrak h_\mu(n)$
defined this way are different from the $\mathfrak a_\mu(n)$
used in \cite{LQW2}.)
Recall further \cite{Na3} that the cohomology degree
for $\mathfrak a_{-k}(\g)$ is $2k -2 +|\g|$.
This implies that the cohomology degree
of $\mathfrak h_\mu(n)$ is exactly $2 \|\mu\|$, and
therefore the cohomology grading and the grading
given by $\|\cdot \|$ are compatible.

\begin{remark}
For each finite subgroup $\G \leq SL_2(\C)$, one
associates the Cartan matrix of the corresponding Dynkin diagram.
The inverse of a Cartan matrix already makes
an appearance in the cup product of the cohomology
class dual to the boundary divisor with itself
in $H^*((\ale)^{[n]})$ for $n \geq 2$. This might help to convince
the reader that the cohomology ring of $(\ale)^{[n]}$
can be very interesting
even though  the cohomology ring of a simply-connected
quasi-projective surface $X$ is trivial.
\end{remark}

\begin{conjecture} (the Constant Conjecture) \label{conj:indep}
Let $X$ be a quasi-projective surface. The structure constants of
the ring $H^*(\Xn)$ with respect to the Heisenberg monomial basis
$(h_\mu(n))$ are independent of $n$. Thus, we can construct a
FH-type ring associated to $X$ which encodes the cohomology rings
$H^*(\Xn)$ for all $n$.
\end{conjecture}

\begin{remark}
A notion of {\em Hilbert ring} $\mathfrak H_X$
was introduced in \cite{LQW2} associated to a {\em projective}
surface $X$ which encodes the cohomology ring structures
of $H^*(\Xn)$ for all $n$. This notion has a counterpart
(which is called the {\em stable ring})
for the orbifold cohomology rings of symmetric products,
and for the class algebras of
wreath products. A central result which
allowed one to introduce these notions is the $n$-independence
of certain coefficients. However, the formulation of
the $n$-independence in those constructions is different
in a fundamental way from the one which we establish in this paper.
The formulation of the Hilbert ring etc is very subtle
as it uses {\em normalized} Heisenberg monomials
for Hilbert schemes and symmetric products
and {\em normalized} conjugacy classes for wreath products,
which are not linearly independent for a fixed $n$.
\end{remark}

\end{document}